# Solving Three Dimensional Maxwell Eigenvalue Problem with Fourteen Bravais Lattices

Tsung-Ming Huang[a], Tiexiang Li[b], Wei-De Li[c], Jia-Wei Lin[d], Wen-Wei Lin[d], Heng Tian[d]

[a]*Department of Mathematics, National Taiwan Normal University, Taipei, 116, Taiwan*
[b]*School of Mathematics, Southeast University, Nanjing 211189, People's Republic of China*
[c]*Department of Mathematics, National Tsing-Hua University, Hsinchu 300, Taiwan*
[d]*Department of Applied Mathematics, National Chiao Tung University, Hsinchu 300, Taiwan*

**Abstract**

Calculation of band structure of three dimensional photonic crystals amounts to solving large-scale Maxwell eigenvalue problems, which are notoriously challenging due to high multiplicity of zero eigenvalue. In this paper, we try to address this problem in such a broad context that band structure of three dimensional isotropic photonic crystals with all 14 Bravais lattices can be efficiently computed in a unified framework. We uncover the delicate machinery behind several key results of our work and on the basis of this new understanding we drastically simplify the derivations, proofs and arguments in our framework. In this work particular effort is made on reformulating the Bloch boundary condition for all 14 Bravais lattices in the redefined orthogonal coordinate system, and establishing eigen-decomposition of discrete partial derivative operators by systematic use of commutativity among them, which has been overlooked previously, and reducing eigen-decomposition of double-curl operator to the canonical form of a $3 \times 3$ complex skew-symmetric matrix under unitary congruence. With the validity of the novel nullspace free method in the broad context, we perform some calculations on one benchmark system to demonstrate the accuracy and efficiency of our algorithm.

## 1. Introduction

The photonic crystal (PhC) is an essential device when light is manipulated in optoelectronics industry. A PhC is a in one-, two- and three-dimensional (1D, 2D, 3D) periodic structure which is composed of different optical media that can purposefully affect electromagnetic wave propagation. This term is coined after Yablonovitch [39] and John [25]'s milestone work in 1987. In recent years, the research about PhC is booming due to the emergence of the topological PhCs (or photonic topological insulator) [33], especially the 3D topological

PhCs. To determine whether a PhC is the topological PhC, band structure calculation is indispensable [28]. To practically know the band structure of a 3D isotropic/anisotropic PhC, we need to first recast the source-free Maxwell equations in frequency domain [37] as follows, with a specific media whose intrinsic properties are described by a 3-by-3 permeability matrix $\mu$ and a permittivity matrix $\varepsilon$, respectively,

$$\nabla \times \boldsymbol{E} = \imath\omega\mu\boldsymbol{H}, \quad \nabla \cdot (\mu\boldsymbol{H}) = 0, \tag{1a}$$

$$\nabla \times \boldsymbol{H} = -\imath\omega\varepsilon\boldsymbol{E}, \quad \nabla \cdot (\varepsilon\boldsymbol{E}) = 0, \tag{1b}$$

where $\omega$ is the frequency, and $\boldsymbol{E}$ and $\boldsymbol{H}$ are the electric and magnetic field, respectively. The famous Bloch theorem [27] requires the solution $\boldsymbol{E}$ and $\boldsymbol{H}$ satisfy the Bloch condition (BC)[34],

$$\boldsymbol{E}(\mathbf{x} + \mathbf{a}_\ell) = \mathbf{e}^{\imath 2\pi\mathbf{k}\cdot\mathbf{a}_\ell}\boldsymbol{E}(\mathbf{x}), \; \boldsymbol{H}(\mathbf{x} + \mathbf{a}_\ell) = \mathbf{e}^{\imath 2\pi\mathbf{k}\cdot\mathbf{a}_\ell}\boldsymbol{H}(\mathbf{x}), \; \ell = 1, 2, 3, \tag{2}$$

where $\mathbf{a}_\ell$ is the lattice translation vectors and $2\pi\mathbf{k}$ is the Bloch wave vector within the first Brillouin zone [23]. For simplicity, $\mu$ is set to the vacuum permeability $\mu_0$ while $\varepsilon$ is assumed to be diagonal throughout this paper.

Given a specific 3D PhC, only certain nonzero $\omega$ can satisfy (1a) (1b) simultaneously. Our ultimate goal is to find a couple of smallest positive eigenvalues of the following Maxwell Eigenvalue Problem (MEP)

$$\begin{bmatrix} & \imath\nabla\times \\ -\imath\nabla\times & \end{bmatrix} \begin{bmatrix} \boldsymbol{E} \\ \boldsymbol{H} \end{bmatrix} = \omega \begin{bmatrix} \varepsilon & \\ & \mu \end{bmatrix} \begin{bmatrix} \boldsymbol{E} \\ \boldsymbol{H} \end{bmatrix}, \tag{3a}$$

$$\nabla \cdot (\varepsilon\boldsymbol{E}) = 0, \quad \nabla \cdot (\mu\boldsymbol{H}) = 0. \tag{3b}$$

To discretize MEP (3), plane-wave expansion method [18, 24, 26, 35], multiple scattering method [16, 36], finite-difference frequency-domain method (FDFD) [9, 10, 15, 20, 21, 38, 40, 41, 42], finite element method [6, 7, 8, 17, 22, 29, 14, 30, 31, 32], to name a few, are available. In the case of diagonal $\varepsilon$ matrix, the Yee's scheme finite-difference scheme [41], originally proposed for time-domain simulation, is particularly attractive. In [20, 21], we have used Yee's scheme [41] for discretization, which results in a generalized eigenvalue problem (GEP). For a 3D PhC, due to divergence-free condition (3b), dimension of the nullspace of the GEP accounts for one third of the total dimension. The presence of the huge nullspace will pose an extraordinary challenge to the desired solutions of the GEP. In fact, no frequency-domain method is immune to this challenge. Besides, even though only smallest few positive eigenvalues are desired, which can be calculated by the Invert-Lanczos method, to solve the corresponding linear system of huge size in each step of the Invert-Lanczos process is another challenge. In [20, 21], we have shown how we resolve these challenges in the case of face-centered cubic (FCC) lattice and simple cubic (SC) lattice.

In this paper, we will solve the MEP (3) for all 14 Bravais lattices along the same lines as [20, 21]. Since the triclinic lattice is the most general one, which can in fact become other 13 Bravais lattices with corresponding constraints imposed, it suffices to consider triclinic lattice only. Several obstacles stand out. For example, since the unit cell of the triclinic lattice is a slanted parallelepiped, it is unclear how to formulate in matrix language the discrete single-curl operator which is compatible with the BC (2), and then it is uncertain whether the advanced nullspace free method in [20] can be applicable in this case. Although it is not uncommon to employ the oblique coordinate system in engineering and physics community, we are not convinced that all our previous inventions can still be applicable in the oblique coordinate system, so we decide to work with the orthogonal coordinate system as before to overcome these obstacles.

This paper is outlined as follows.

- In Sec. 2 an orthogonal coordinate system with which we actually work are built from the oblique coordinate system generated by $\mathbf{a}_1, \mathbf{a}_2, \mathbf{a}_3$.
- In Sec. 3 we reformulate the BC (2) within the cubic working cell.
- In Sec. 4 we discretize $\nabla \times E$ and $\nabla \times H$ into matrix-vector products,

  $$\mathcal{C}E = \begin{bmatrix} 0 & -C_3 & C_2 \\ C_3 & 0 & -C_1 \\ -C_2 & C_1 & 0 \end{bmatrix} \begin{bmatrix} E_1(:) \\ E_2(:) \\ E_3(:) \end{bmatrix}, \mathcal{C}^* H = \begin{bmatrix} 0 & C_3^* & -C_2^* \\ -C_3^* & 0 & C_1^* \\ C_2^* & -C_1^* & 0 \end{bmatrix} \begin{bmatrix} H_1(:) \\ H_2(:) \\ H_3(:) \end{bmatrix},$$

  respectively, where $C_2, C_3$ are quite complicated in most Bravais lattices due to the reformulated BC. Moreover, we reduce the MEP (3) into a GEP: $\mathcal{A}E = \lambda \mathcal{B}E,\ \mathcal{A} = \mathcal{C}^*\mathcal{C},\ \lambda = \mu_0\omega^2$.
- In Sec. 5 we prove that $C_1, C_2, C_3$ commute with each othe, and obtain eigen-decomposition of them analytically: $C_1 T = T\Lambda_1,\ C_2 T = T\Lambda_2,\ C_3 T = T\Lambda_3$, where $\Lambda_1, \Lambda_2, \Lambda_3$ are diagonal matrices.
- In Sec. 6 we analytically identify the orthonormal basis of nullspace and range space of $\mathcal{C}$ for any given Bravais lattice, and set up eigen-decomposition of the discrete double-curl operator $\mathcal{A}$ *i.e.*, $\mathcal{A} = \mathcal{Q}_r \Lambda_r \mathcal{Q}_r^*,\ Q_r^* Q_r = I_{2n}$, in light of the canonical form of a $3 \times 3$ complex skew-symmetric matrix under unitary congruence: $L = \overline{\widetilde{V}} \begin{bmatrix} 0 & -\beta \\ \beta & 0 \end{bmatrix} \widetilde{V}^*,\ L^\top = -L \in \mathbb{C}^{3\times 3},\ \widetilde{V}^*\widetilde{V} = I_2$.
- In Sec. 7 by eliminating the considerable nullspace of $\mathcal{A}$, we transform the GEP into a nullspace free standard eigenvalue problem (NFSEP): $\mathcal{A}_r \widetilde{E} = \lambda \widetilde{E},\ \mathcal{A}_r = \mathcal{A}_r^* > 0$. For the sake of self-containedness of this article, fast eigensolver for NFSEP is reviewed.
- In Sec. 8, some numerical results are presented to demonstrate the accuracy and efficiency of our method.
- Finally, we conclude our present work in Sec. 9.

Here we briefly introduce some notations commonly used in this work. $A^*, A^\top$ denote conjugate transpose and transpose of matrix $A$, respectively, and $\bar{\cdot}$ denotes complex conjugate. $||\cdot||_2$ denotes the Euclidean norm. A 3D vector is marked in bold and is equivalent to its Cartesian coordinate representation. $\otimes$ denotes the Kronecker product. $A \oplus B$ means the direct sum of matrices $A, B$. $\delta_{\ell',\ell}$ is the Kronecker delta function, *i.e.*, $\delta_{\ell',\ell} = 1$ if $\ell' = \ell$ and $\delta_{\ell',\ell} = 0$ otherwise. $e_\ell = [\delta_{1,\ell},\ \delta_{2,\ell}, \cdots \delta_{n,\ell}]^\top$ is the standard unit vector in $\mathbb{R}^n$. We define $\xi(\theta) := \exp(\imath 2\pi\theta)$. $\square$ABCD refers to rectangular ABCD. For convenience we will employ Matlab[3] language with little explanation. For example, 'floor' denotes the function of rounding to the nearest integer towards $-\infty$. Let $\mathrm{vec}(X)$ denote the vectorization operation of a matrix $X$ of any size, *i.e.*, $\mathrm{vec}(X) = X(:)$.

## 2. Lattice translation vectors, physical cell and working cell

A crystal structure can be regarded as a lattice structure plus a basis. At present, millions of crystals are known, and each crystal has a different nature. Fortunately, there are only 7 lattice systems and 14 Bravais lattices in 3D Euclidean space [1]. The so-called primitive cell is a fundamental domain under the translational symmetry, and contains just one lattice point[4]. In fact a 3D primitive cell is a slanted parallelepiped formed by lattice translation vectors $\mathbf{a}_1, \mathbf{a}_2$ and $\mathbf{a}_3$, as illustrated in Figure 1(a). As mentioned above, in triclinic lattice there is no restriction on the length of $\mathbf{a}_1, \mathbf{a}_2, \mathbf{a}_3$ nor on the angle between any two of them, if we are able to solve the MEP (3) in triclinic lattice, we can also cope with other lattices in almost the same manner. Therefore we will focus on the triclinic lattice in the main body of this work. For convenience, we dub the primitive cell of triclinic lattice as 3D physical cell. In that it is inconvenient to discretize MEP (3) in the 3D physical cell using finite difference, we need to redefine a cuboid primitive cell generated by new vectors $\mathbf{a}, \mathbf{b}, \mathbf{c}$ which are orthogonal basis of $\mathbf{a}_1, \mathbf{a}_2, \mathbf{a}_3$. Specifically, we first sort $\mathbf{a}_1, \mathbf{a}_2, \mathbf{a}_3$ in terms of length in descending order, resulting in $\tilde{\mathbf{a}}_1, \tilde{\mathbf{a}}_2, \tilde{\mathbf{a}}_3$. Then we require

$$\mathbf{a} = \tilde{\mathbf{a}}_1, \quad \mathbf{c} /\!/ \tilde{\mathbf{a}}_1 \times \tilde{\mathbf{a}}_2, \quad \mathbf{b} /\!/ \mathbf{c} \times \mathbf{a}, \tag{4}$$

$$\mathbf{b} \times \mathbf{a} = \tilde{\mathbf{a}}_2 \times \tilde{\mathbf{a}}_1, \quad \mathbf{c} \cdot \mathbf{a} \times \mathbf{b} = \tilde{\mathbf{a}}_3 \cdot \tilde{\mathbf{a}}_1 \times \tilde{\mathbf{a}}_2. \tag{5}$$

It is easy to see in Figure 1(a) that $\mathbf{a}, \mathbf{b}, \mathbf{c}$ defined in this way are unique. Let $a, b, c$ be the lengths of $\mathbf{a}, \mathbf{b}, \mathbf{c}$, respectively. Identifying $\mathbf{a}/a, \mathbf{b}/b, \mathbf{c}/c$ as unit vectors of the usual $x$-,$y$-,$z$-axis, vectors $\tilde{\mathbf{a}}_1, \tilde{\mathbf{a}}_2, \tilde{\mathbf{a}}_3$ can be rewritten as

$$\left[\tilde{\mathbf{a}}_1, \tilde{\mathbf{a}}_2, \tilde{\mathbf{a}}_3\right] = \begin{bmatrix} a_1 & a_2 \cos\phi_3 & a_3 \cos\phi_2 \\ 0 & a_2 \sin\phi_3 & a_3 \ell_2 \\ 0 & 0 & a_3 \ell_3 \end{bmatrix}, \tag{6}$$

where $a_j$ is the length of $\tilde{\mathbf{a}}_j$ and $\phi_j$ is angle between $\tilde{\mathbf{a}}_i$ and $\tilde{\mathbf{a}}_k$, $i, j, k = 1, 2, 3, i \neq j \neq k$, and

$$\ell_2 = \left(\cos\phi_1 - \cos\phi_3 \cos\phi_2\right) / \sin\phi_3, \tag{7a}$$

$$\ell_3 = \sqrt{\sin^2\phi_2 - \ell_2^2}, \tag{7b}$$

with $a_2 \sin\phi_3 > |\ell_2|$. This result, which can also be found in [2], is illustrated in Figure 1(b). In passing, we have

$$a = a_1, \quad b = a_2 \sin\phi_3, \quad c = a_3\ell_3. \tag{8}$$

We tabulate the coordinate representation (6) of $\tilde{\mathbf{a}}_1, \tilde{\mathbf{a}}_2, \tilde{\mathbf{a}}_3$ of all 7 lattice systems in Appendix A. We will solve MEP (3) mainly in this cuboid primitive cell which is dubbed as 3D working cell. To convey the basic techniques and methods in our framework of modeling of 3D PhCs, we just work on one specific case where $\phi_3 < \pi/2,\ \phi_2 < \pi/2,\ \ell_2 > 0$ in the main body of this article, and defer the discussion of other possible combination of $\phi_3,\ \phi_2,\ \ell_2$ to the Appendix.

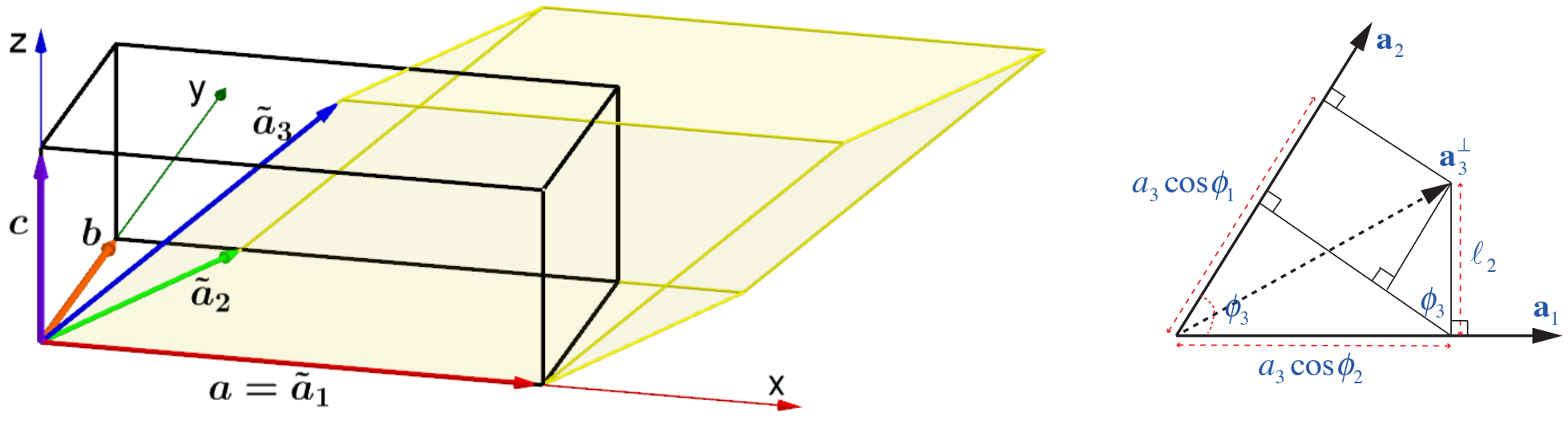

(a) triclinic lattice

(b) $\mathbf{a}_3^{\perp} = \mathbf{a}_3 - \mathbf{c}, \mathbf{c} \perp \mathbf{a}_3^{\perp}$

Figure 1: triclinic lattice and its projective view

## 3. BC (2) within the working cell

Hereafter, for simplicity, we assume $\tilde{\mathbf{a}}_1, \tilde{\mathbf{a}}_2, \tilde{\mathbf{a}}_3$ are just $\mathbf{a}_1, \mathbf{a}_2, \mathbf{a}_3$. Viewed in the associated oblique coordinate system spanned by $\mathbf{a}_1, \mathbf{a}_2, \mathbf{a}_3$, BC (2) is very clear, and is naturally compatible with periodicity of a PhC along $\mathbf{a}_1, \mathbf{a}_2, \mathbf{a}_3$. However, in the working cell or the orthogonal coordinate system with $x$-,$y$-,$z$-axis, formulation of BC needs some effort.

Given $\mathbf{v} \in \mathbb{R}^3$, the translation operator $\mathcal{T}_{\mathbf{v}}$ is defined as $\mathcal{T}_{\mathbf{v}}(\mathbf{x}) := \mathbf{x} + \mathbf{v}, \forall \mathbf{x} \in \mathbb{R}^3$. Clearly, $\mathcal{T}_{\mathbf{v}_1+\mathbf{v}_2} = \mathcal{T}_{\mathbf{v}_1}\mathcal{T}_{\mathbf{v}_2} = \mathcal{T}_{\mathbf{v}_2}\mathcal{T}_{\mathbf{v}_1}$.

In fact, the 3D working cell is the set $\mathbb{D} = [0, a) \times [0, b) \times [0, c) \subset \mathbb{R}^3$, *i.e.*, a cuboid of lengths $a, b, c$. Since $\mathbf{a}_1 = \mathbf{a}$, BC along $x$-axis is just the same. Specifically, given $\mathbf{x} = (x, y, z) \in \mathbb{D}$,

$$\boldsymbol{E}(\mathbf{x}) = \xi(\mathbf{k} \cdot (\mathbf{x} - \mathcal{T}_{-\mathbf{a}}(\mathbf{x})))\boldsymbol{E}(\mathcal{T}_{-\mathbf{a}}(\mathbf{x})). \tag{9}$$

However, the relation between $\boldsymbol{E}(\mathbf{x})$ and $\boldsymbol{E}(\mathcal{T}_{-\mathbf{b}}(\mathbf{x}))$ or $\boldsymbol{E}(\mathcal{T}_{-\mathbf{c}}(\mathbf{x}))$ can not resemble (9) naively. Fortunately, for derivations in Sec. 4, we only need to know $\boldsymbol{E}(\mathcal{T}_{-\mathbf{c}}((x, y, c)))$ with $(x, y, c) \in \mathbb{D}$ and $\boldsymbol{E}(\mathcal{T}_{-\mathbf{b}}((x, b, z))))$ with $(x, b, z) \in \mathbb{D}$.

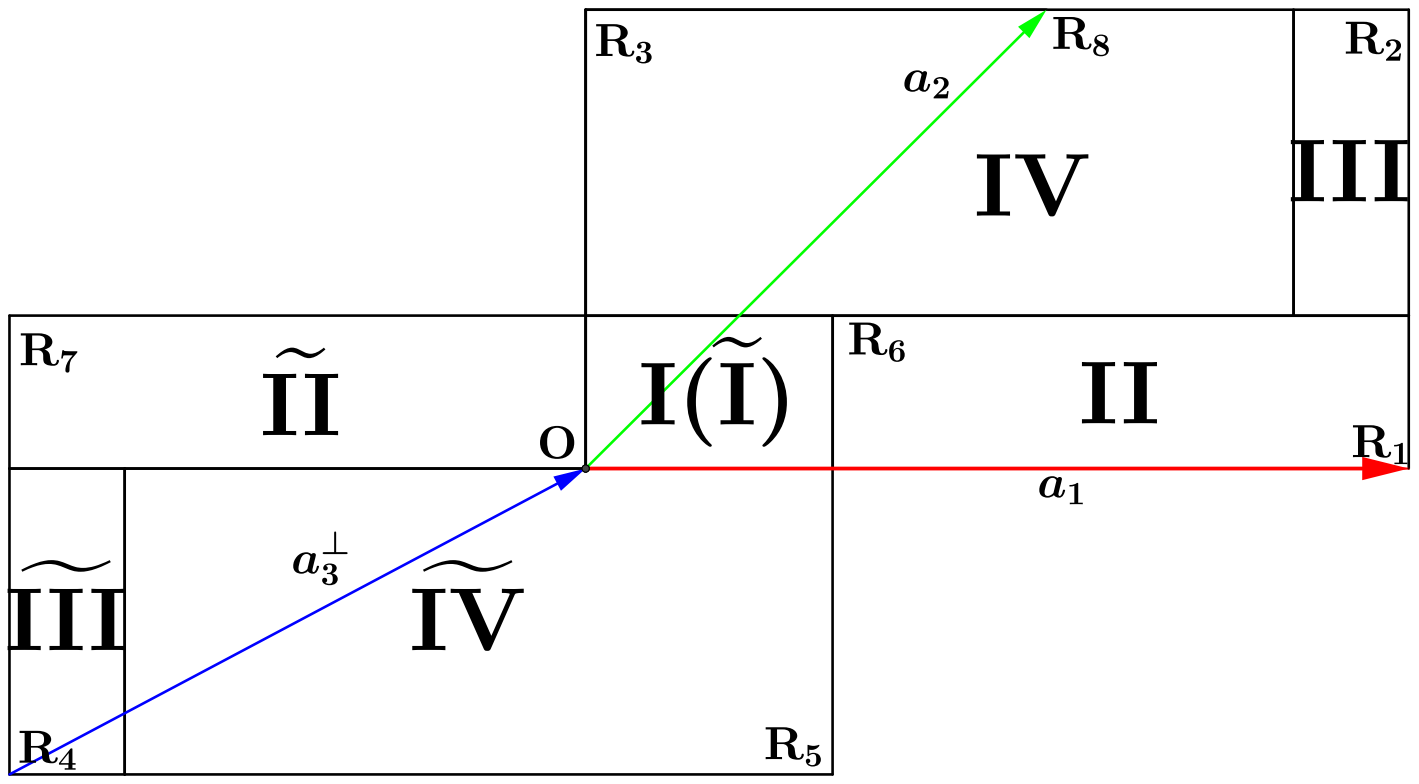

Figure 2: illustration of $(x, y, 0) \in \mathbb{D}$,*i.e.*, the bottom surface, and $\mathcal{T}_{-\mathbf{a}_3}((x, y, c))$

Given $\mathbf{x} = (x, y, 0) \in \mathbb{D}$, we just think of $(x_2, y_2, 0)$ as image of $\mathbf{x} + c$, *i.e.*, a point of the top surface of $\mathbb{D}$, under $\mathcal{T}_{-\mathbf{a}_3}$, then BC along $z$-direction could be

$$\begin{aligned} \boldsymbol{E}(\mathcal{T}_{-\mathbf{c}}((x, y, c))) = \boldsymbol{E}(((x, y, 0)) &= \xi\left(\mathbf{k} \cdot ((x, y, 0)) - (x_2, y_2, 0)\right) \boldsymbol{E}((x_2, y_2, 0)) \\ &= \xi\left(\mathbf{k} \cdot ((x, y, 0)) - \mathcal{T}_{-\mathbf{a}_3}((x, y, c))\right) \boldsymbol{E}(\mathcal{T}_{-\mathbf{a}_3}((x, y, c))), \end{aligned} \tag{10}$$

with $(x, y, 0) - \mathcal{T}_{-\mathbf{a}_3}((x, y, c))$ being integer multiples of $\mathbf{a}_1, \mathbf{a}_2$.

In Figure 2, $\square \mathrm{OR}_1\mathrm{R}_2\mathrm{R}_3$ is the bottom surface of $\mathbb{D}$, while $\square \mathrm{R}_4\mathrm{R}_5\mathrm{R}_6\mathrm{R}_7$ is image of the top surface of $\mathbb{D}$ under $\mathcal{T}_{-\mathbf{a}_3}$ and overlaps with patch I of the former. In short, there should be four patches within $\square \mathrm{OR}_1\mathrm{R}_2\mathrm{R}_3$, namely, $\mathrm{I}, \mathrm{II}, \mathrm{III}, \mathrm{IV}$, and each patch, equipped with different linear mapping, namely, $\mathcal{T}_0, \mathcal{T}_{-\mathbf{a}_1}, \mathcal{T}_{-\mathbf{a}_1-\mathbf{a}_2}, \mathcal{T}_{-\mathbf{a}_2}$, respectively, is mapped to four patches, namely, $\widetilde{\mathrm{I}}, \widetilde{\mathrm{II}}, \widetilde{\mathrm{III}}, \widetilde{\mathrm{IV}}$, respectively, within $\square \mathrm{R}_4\mathrm{R}_5\mathrm{R}_6\mathrm{R}_7$. Then we can establish the correct BC within the bottom surface of $\mathbb{D}$, which specifies $x_2, y_2$ in (10). Letting $\mathbf{x} = (x, y, 0) \in \mathbb{D}$,

$$\boldsymbol{E}(\mathbf{x}) = \begin{cases} \boldsymbol{E}(\mathbf{x}), & \text{if } \mathbf{x} \in \mathrm{I} \\ \xi(\mathbf{k} \cdot \mathbf{a}_1)\boldsymbol{E}(\mathbf{x} - \mathbf{a}_1), & \text{if } \mathbf{x} \in \mathrm{II} \\ \xi(\mathbf{k} \cdot (\mathbf{a}_1 + \mathbf{a}_2))\boldsymbol{E}((\mathbf{x} - \mathbf{a}_1 - \mathbf{a}_2)), & \text{if } \mathbf{x} \in \mathrm{III} \\ \xi(\mathbf{k} \cdot \mathbf{a}_2)\boldsymbol{E}(\mathbf{x} - \mathbf{a}_2), & \text{if } \mathbf{x} \in \mathrm{IV}. \end{cases} \tag{11}$$

In passing, considering that $\boldsymbol{E}(\mathcal{T}_{\mathbf{a}_3}(\mathbf{x})) = \xi(\mathbf{k} \cdot \mathbf{a}_3)\boldsymbol{E}(\mathbf{x})$, we can of course add $\mathbf{a}_3$ to the argument of $\boldsymbol{E}$ on the right hand side of (11) with updated prefactor.

We refer the reader to Appendix B to see how Figure 2 as well as the associated mapping is obtained. Depending on different combination of $\phi_3$, $\phi_2$, $\ell_2$, BC (11) could be different. Results in other cases will be listed in Appendix.

As for $\boldsymbol{E}(\mathcal{T}_{-\mathbf{b}}((x, b, z))))$ with $(x, b, z) \in \mathbb{D}$, since $z$-axis is independent of $x$-,$y$-axis, we can just let $z = 0$ here for simplicity. Letting $\mathbf{x} = (x, b, 0) \in \mathbb{D}$, we have BC along $y$-direction for different patches of $\mathrm{R}_3\mathrm{R}_2$ shown in Figure 2:

$$\boldsymbol{E}(\mathbf{x}) = \begin{cases} \xi(\mathbf{k} \cdot \mathbf{a}_2)\boldsymbol{E}(\mathcal{T}_{-\mathbf{a}_2}(\mathbf{x})), & \text{if } \mathbf{x} \in \mathrm{R}_8\mathrm{R}_2 \\ \xi(\mathbf{k} \cdot (\mathbf{a}_2 - \mathbf{a}_1))\boldsymbol{E}(\mathcal{T}_{\mathbf{a}_1-\mathbf{a}_2}(\mathbf{x})), & \text{if } \mathbf{x} \in \mathrm{R}_3\mathrm{R}_8. \end{cases} \tag{12}$$

## 4. Matrix Representation of the Discretized Single-Curl

As mentioned above, partial derivative of the trivariate function $\boldsymbol{E}(\mathbf{x}), \boldsymbol{H}(\mathbf{x})$ in (1a) and (1b) which satisfies 3D BC (2) will be approximated by finite difference of this function.

First of all, finite difference approximation needs grid, usually the uniform grid, in the simulation domain. Given $n_1, n_2, n_3 \in \mathbb{N}$, we can have a uniform grid along $x$-,$y$-,$z$-axis of our 3D working cell $\mathbb{D}$, defined in Sec. 3, respectively, with constant grid spacing

$$\delta_x = a/n_1, \quad \delta_y = b/n_2, \quad \delta_z = c/n_3,$$

respectively. In general, each component of the vector-valued function $\boldsymbol{E}(\mathbf{x}) = [E_1(\mathbf{x}), E_2(\mathbf{x}), E_3(\mathbf{x})]^\top$ could be sampled at different points. Hence we assume that $E_\ell(\mathbf{x})$ is sampled at

$$\mathbf{x}_\ell(i,j,k) = \mathbf{x}_\ell(0,0,0) + (i\delta_x, j\delta_y, k\delta_z), \tag{13}$$

where $\mathbf{x}_\ell(0,0,0)$ will be specified later in this section and $\ell = 1,2,3,\ i = 0,1,\ldots,n_1-1,\ j = 0,1,\ldots,n_2-1,\ k = 0,1,\ldots,n_3-1$. Unless otherwise stated, in this section $i,j,k$ always take on these values. In passing, this definition (13) of $\mathbf{x}_\ell(i,j,k)$ holds even when $i,j,k \in \mathbb{Z}$.

Given $\ell$, the three-way array $E_\ell(\mathbf{x}_\ell(:,:,:))$ is arranged in column-major order, *i.e.*, the first index varies fastest while the last index varies slowest. For convenience we store $E_\ell(\mathbf{x}_\ell(:,:,:))$ in a column vector

$$E = [E_1(:); E_2(:); E_3(:)].$$

Let's deal with single-curl $\nabla\times$ in (1a) first, without worrying about $\partial_x E_1$ *etc.* at the moment. Below we will refer to quantities in (6) and (8) frequently.

**Part I. Discrete $\partial_x E_\ell$.** Since BC (9) is very similar to 1D case, using matrix language, we recast

$$\frac{E_\ell(\mathbf{x}_\ell(i+1,j,k)) - E_\ell(\mathbf{x}_\ell(i,j,k))}{\delta_x}, \quad \ell = 2,3, \tag{14}$$

into $C_1 E_\ell(:)$, where

$$C_1 = I_{n_3} \otimes I_{n_2} \otimes \frac{K_1 - I_{n_1}}{\delta_x} \in \mathbb{C}^{n\times n}, \tag{15}$$

$$K_1 = \begin{bmatrix} 0 & I_{n_1-1} \\ \xi(\mathbf{k}\cdot\mathbf{a}_1) & 0 \end{bmatrix} \in \mathbb{C}^{n_1\times n_1}. \tag{16}$$

**Part II. Discrete $\partial_y E_\ell$.** BC (12) holds for continuous $\mathbf{x}$, however if we want to recast

$$\frac{E_\ell(\mathbf{x}_\ell(i,j+1,k)) - E_\ell(\mathbf{x}_\ell(i,j,k))}{\delta_y}, \quad \ell = 1,3, \tag{17}$$

into matrix-vector product, we need the discretized version of BC (12).

Although in Figure 3 we have in principle $\mathrm{R}_8 \equiv \mathrm{O} \mod \mathbf{a}_2$, it is very rare that $\mathrm{R}_8$ coincides exactly with any of the grid point in a given uniform grid in $\mathrm{R}_3\mathrm{R}_2$. As an expediency to resolve this mismatching, we stipulate that the rightmost grid point within $\mathrm{R}_3\mathrm{R}_8$ is the substitute of $\mathrm{R}_8$. Putting it differently, when $\phi_3 < \pi/2$, since the number of grid points in $\mathrm{R}_3\mathrm{R}_8$ is $m_1$, where

$$m_1 = \text{floor}\left((\mathbf{a}_2 \cdot \mathbf{a}_1/a_1)/\delta_x\right) = \text{floor}\left((a_2 \cos\phi_3)/\delta_x\right), \tag{18}$$

then

$$\mathbf{x}_\ell(m_1, n_2, k) \equiv \mathbf{x}_\ell(0, 0, k) \mod \mathbf{a}_2,$$

holds at the discrete level by force.

$E_\ell(\mathbf{x}_\ell(:, n_2, k))$, a column vector of length $n_1$, is partitioned into 2 blocks, due to two cases in (12). Then the discretized BC (12) is

$$E_\ell(\mathbf{x}_\ell(:, n_2, k)) = \xi(\mathbf{k} \cdot \mathbf{a}_2) J_2 E_\ell(\mathbf{x}_\ell(:, 0, k)), \tag{19}$$

$$J_2 = \begin{bmatrix} 0 & \xi(-\mathbf{k} \cdot \mathbf{a}_1) I_{m_1} \\ I_{n_1 - m_1} & 0 \end{bmatrix}. \tag{20}$$

Finally, (17) is recast into $C_2 E_\ell(:)$, where

$$C_2 = I_{n_3} \otimes \frac{K_2 - I_{n_1 n_2}}{\delta_y} \in \mathbb{C}^{n \times n}, \tag{21}$$

$$K_2 = \begin{bmatrix} 0 & I_{n_2 - 1} \otimes I_{n_1} \\ \xi(\mathbf{k} \cdot \mathbf{a}_2) J_2 & 0 \end{bmatrix} \in \mathbb{C}^{(n_1 n_2) \times (n_1 n_2)}. \tag{22}$$

In the case where $\phi_3 > \pi/2$, the expression for $m_1, J_2$ can be found in Appendix.

**Part III. Discrete $\partial_z E_\ell$.** If we want to recast

$$\frac{E_\ell(\mathbf{x}_\ell(i, j, k+1)) - E_\ell(\mathbf{x}_\ell(i, j, k))}{\delta_z}, \quad \ell = 1, 2, \tag{23}$$

into a matrix-vector product, we need to know how $E_\ell(\mathbf{x}_\ell(:, :, n_3))$ is related to $E_\ell(\mathbf{x}_\ell(:, :, 0))$ from BC (11).

We have the following observations about Figure 3,

- length of $\mathrm{R}_9\mathrm{R}_6$ is $a_1 - a_3 \cos\phi_2$, while length of $\mathrm{R}_9\widehat{\mathrm{R}}_5$ is $a_3 \cos\phi_2 - a_2 \cos\phi_3$.
- length of $\mathrm{R}_3\mathrm{R}_9$ is $a_3 \ell_2$, while length of $\mathrm{R}_9\mathrm{O}$ is $a_2 \sin\phi_3 - a_3 \ell_2$.

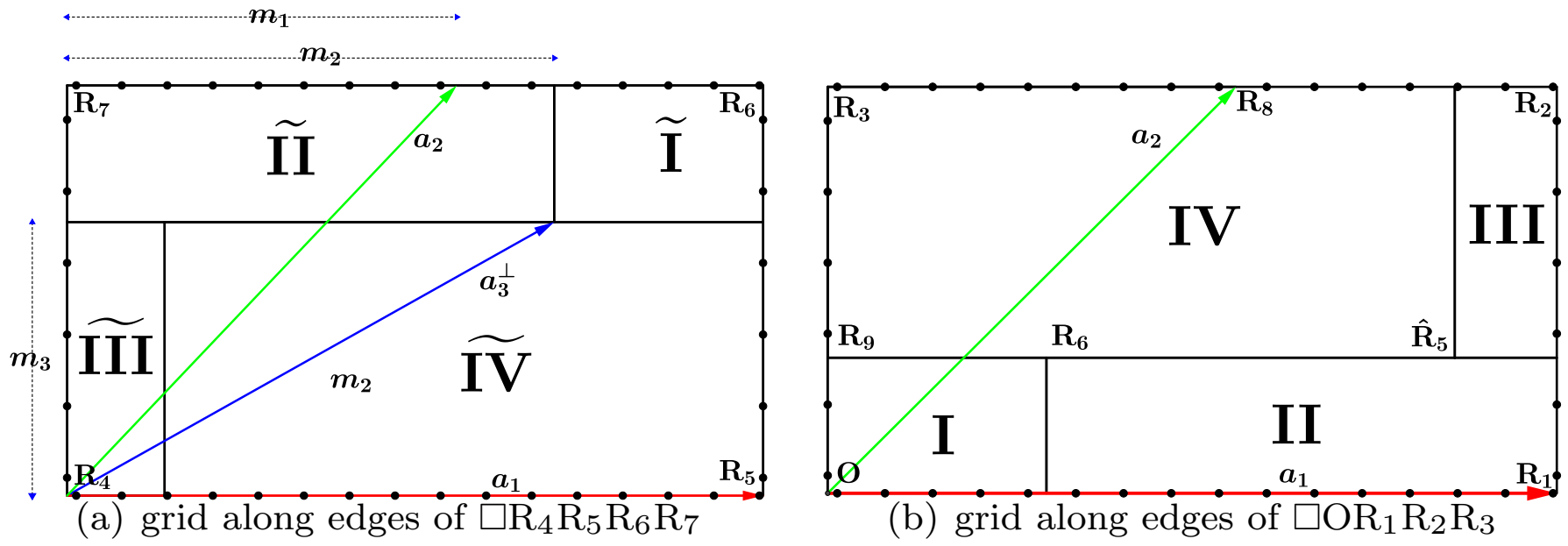

(a) grid along edges of $\square R_4R_5R_6R_7$ (b) grid along edges of $\square OR_1R_2R_3$

Figure 3: illustration of uniform grid within top and bottom surface of $\mathbb{D}$.

Again, it is very rare that vertices of any patch in Figure 3 coincide exactly with any of the grid point for a given uniform mesh in $\square OR_1R_2R_3$. Define

$$m_2 = \text{floor}((\mathbf{a}_3 \cdot \mathbf{a}_1/a_1)/\delta_x) = \text{floor}((a_3 \cos\phi_2)/\delta_x), \tag{24}$$

$$m_3 = \text{floor}((\mathbf{a}_3 \cdot \mathbf{b}/b)/\delta_y) = \text{floor}(a_3\ell_2/\delta_y), \tag{25}$$

$$m_4 = m_2 - m_1, \tag{26}$$

then along $x$-direction, $R_9R_6$ contains $n_1 - m_2$ grid points, $R_9\widehat{R}_5$ contains $n_1 + m_1 - m_2$ grid points, while along $y$-direction, $R_3R_9$ contains $m_3$ grid points, the edge $R_9O$ contains $n_2 - m_3$ grid points.

Matrices $E_\ell(\mathcal{T}_{-\mathbf{a}_3}(\mathbf{x}_\ell(:,:,n_3)))$ and $E_\ell(\mathbf{x}_\ell(:,:,0))$ are partitioned into 4 blocks,

$$E_\ell(\mathbf{x}_\ell(:,:,0)) = \begin{bmatrix} E_{\text{I}} & E_{\text{IV}} \\ E_{\text{II}} & E_{\text{III}} \end{bmatrix} \in \mathbb{C}^{n_1\times n_2}, \quad E_\ell(\mathcal{T}_{-\mathbf{a}_3}(\mathbf{x}_\ell(:,:,n_3))) = \begin{bmatrix} E_{\widetilde{\text{III}}} & E_{\widetilde{\text{II}}} \\ E_{\widetilde{\text{IV}}} & E_{\widetilde{\text{I}}} \end{bmatrix} \in \mathbb{C}^{n_1\times n_2},$$

in accordance with Figure 3, size of each block of which becomes transparent in (27),(28),(29) below. Then the discretized version of BC (11) is as follows:

$$\begin{bmatrix} E_{\widetilde{\text{II}}} \\ E_{\widetilde{\text{I}}} \end{bmatrix} = \begin{bmatrix} 0 & \xi(-\mathbf{k}\cdot\mathbf{a}_1)I_{m_2} \\ I_{n_1-m_2} & 0 \end{bmatrix} \begin{bmatrix} E_{\text{I}} \\ E_{\text{II}} \end{bmatrix} I_{n_2-m_3}, \tag{27}$$

$$\begin{bmatrix} E_{\widetilde{\text{III}}} \\ E_{\widetilde{\text{IV}}} \end{bmatrix} = \begin{bmatrix} 0 & \xi(-\mathbf{k}\cdot\mathbf{a}_1)I_{m_4} \\ I_{n_1-m_4} & 0 \end{bmatrix} \begin{bmatrix} E_{\text{IV}} \\ E_{\text{III}} \end{bmatrix} \xi(-\mathbf{k}\cdot\mathbf{a}_2)I_{m_3}, \tag{28}$$

$$\begin{bmatrix} E_{\text{IV}} & E_{\text{I}} \\ E_{\text{III}} & E_{\text{II}} \end{bmatrix} = I_{n_1} \begin{bmatrix} E_{\text{I}} & E_{\text{IV}} \\ E_{\text{II}} & E_{\text{III}} \end{bmatrix} \begin{bmatrix} 0 & I_{n_2-m_3} \\ I_{m_3} & 0 \end{bmatrix}. \tag{29}$$

Actually $\text{vec}(E_\ell(\mathbf{x}_\ell(:,:,0)))$ can be seen as vertical concatenation of $\text{vec}\left([E_{\text{I}}; E_{\text{II}}]\right)$ and $\text{vec}\left([E_{\text{IV}}; E_{\text{III}}]\right)$, while $\text{vec}(E_\ell(\mathcal{T}_{-\mathbf{a}_3}(\mathbf{x}_\ell(:,:,n_3))))$ can be viewed as vertical concatenation of $\text{vec}\left([E_{\widetilde{\text{III}}}; E_{\widetilde{\text{IV}}}]\right)$ and $\text{vec}\left([E_{\widetilde{\text{II}}}; E_{\widetilde{\text{I}}}]\right)$. It is easy to see that

$$\left(Z^\top \otimes Y\right)\text{vec}(X) = \text{vec}(YXZ), \tag{30}$$

where $X, Y, Z$ are matrices with any compatible size.

Finally, with (27),(28),(29),(30), we can recast (23) into $C_3 E_\ell(:)$, where

$$
C_3 = \frac{K_3 - I_n}{\delta_z}, \quad K_3 = \begin{bmatrix} 0 & I_{n_3-1} \otimes I_{n_2} \otimes I_{n_1} \\ \xi(\mathbf{k} \cdot \mathbf{a}_3) J_3 & 0 \end{bmatrix} \in \mathbb{C}^{n \times n}, \tag{31}
$$

$$
\begin{aligned}
J_3 &= \left( \begin{bmatrix} 0 & I_{n_2-m_3} \\ I_{m_3} & 0 \end{bmatrix}^\top \otimes I_{n_1} \right) \times \\
&\quad \begin{bmatrix} I_{n_2-m_3} \otimes \begin{bmatrix} 0 & \xi(-\mathbf{k} \cdot \mathbf{a}_1) I_{m_2} \\ I_{n_1-m_2} & 0 \end{bmatrix} & \\ & \xi(-\mathbf{k} \cdot \mathbf{a}_2) I_{m_3} \otimes \begin{bmatrix} 0 & \xi(-\mathbf{k} \cdot \mathbf{a}_1) I_{m_4} \\ I_{n_1-m_4} & 0 \end{bmatrix} \end{bmatrix} \\
&= \begin{bmatrix} & \xi(-\mathbf{k} \cdot \mathbf{a}_2) I_{m_3} \otimes \begin{bmatrix} 0 & \xi(-\mathbf{k} \cdot \mathbf{a}_1) I_{m_4} \\ I_{n_1-m_4} & 0 \end{bmatrix} \\ I_{n_2-m_3} \otimes \begin{bmatrix} 0 & \xi(-\mathbf{k} \cdot \mathbf{a}_1) I_{m_2} \\ I_{n_1-m_2} & 0 \end{bmatrix} & \end{bmatrix}.
\end{aligned} \tag{32}
$$

Depending on different combination of $\phi_3,\ \phi_2,\ \ell_2$, matrix $J_3$ could be different. We put the result for $J_3$ in other cases in Appendix D.

**Part IV. Discrete $\partial_x H_\ell, \partial_y H_\ell, \partial_z H_\ell$.** In order to preserve the Hermiticity of the operator on the left hand side of MEP (3) at the discrete level, single-curl operator in (1b) should be discretized slightly differently. We will not detail the derivations, but just present the results. Specifically, the discretized version of BC (9),(12) and (11) can be immediately written down verbatim in terms of $\boldsymbol{H}(\mathbf{x})$ in place of $\boldsymbol{E}(\mathbf{x})$, then we can recast

$$
\frac{H_\ell(\mathbf{x}_\ell(i,j,k)) - H_\ell(\mathbf{x}_\ell(i-1,j,k))}{\delta_x}, \quad \ell = 2, 3, \tag{33}
$$

$$
\frac{H_\ell(\mathbf{x}_\ell(i,j,k)) - H_\ell(\mathbf{x}_\ell(i,j-1,k))}{\delta_y}, \quad \ell = 1, 3, \tag{34}
$$

$$
\frac{H_\ell(\mathbf{x}_\ell(i,j,k)) - H_\ell(\mathbf{x}_\ell(i,j,k-1))}{\delta_z}, \quad \ell = 1, 2, \tag{35}
$$

into $-C_1^* H_\ell(:),\ -C_2^* H_\ell(:)$ and $-C_3^* H_\ell(:)$, respectively.

**Part V. Yee's scheme.** To return to the famous Yee's scheme for $\boldsymbol{E}(\mathbf{x})$, $\mathbf{x}_\ell(0,0,0)$ in (13) is set to

$$
\mathbf{x}_1(0,0,0) = (\delta_x/2, 0, 0),\ \mathbf{x}_2(0,0,0) = (0, \delta_y/2, 0),\ \mathbf{x}_3(0,0,0) = (0, 0, \delta_z/2).
$$

In addition, since $\varepsilon(\mathbf{x})$ is assumed to be diagonal, then we can define the following positive diagonal matrix $\mathcal{B}$,

$$
\mathcal{B} = \operatorname{diag}([\operatorname{vec}(\varepsilon(\mathbf{x}_1(:,:,:))); \operatorname{vec}(\varepsilon(\mathbf{x}_2(:,:,:))); \operatorname{vec}(\varepsilon(\mathbf{x}_3(:,:,:)))]),
$$

in which $\mathbf{x}_\ell(i,j,k)$ coincide the grid points where $\boldsymbol{E}(\mathbf{x})$ is sampled.

To return to the Yee's scheme for $\boldsymbol{H}(\mathbf{x})$, $\mathbf{x}_\ell(0,0,0)$ in (13) is set to

$$\mathbf{x}_\ell(0,0,0) = ((1-\delta_{1,\ell})\delta_x/2, (1-\delta_{2,\ell})\delta_y/2, (1-\delta_{3,\ell})\delta_z/2).$$

With Yee's scheme $\mathbf{x}_\ell(i,j,k)$ for $\boldsymbol{E}(\mathbf{x})$ and $\boldsymbol{H}(\mathbf{x})$ specified above, using (14),(17),(23) and (33),(34),(35), it can be proved that the divergence free condition (3b) is automatically satisfied. This is where the superiority of Yee's scheme lies.

**Part VI. Discrete MEP** (3)**.** At last, the discretization of (3) is reduced into the following GEP with size halved,

$$\mathcal{A}E = \lambda \mathcal{B}E, \quad \lambda = \mu_0\omega^2, \quad \mathcal{A} = \mathcal{C}^*\mathcal{C}, \tag{36}$$

$$\mathcal{C} = \begin{bmatrix} 0 & -C_3 & C_2 \\ C_3 & 0 & -C_1 \\ -C_2 & C_1 & 0 \end{bmatrix}. \tag{37}$$

## 5. Eigen-decomposition of partial derivative operators

It is known in Ref. [20, 21] for the FCC lattice without eigen-decomposition of $K_1, K_2, K_3$, it is unlikely to obtain eigen-decomposition of $\mathcal{A}$ in (36) analytically, let alone the nullspace of $\mathcal{A}$. This is also the case for the triclinic lattice and other lattices. The derivation in Ref. [20, 21] could be applied to our present problem with necessary modification, but it turns out the whole process is very complicated and error-prone. In addition, the derivation there can not explain why we have Kronecker product decomposition of $K_2$'s and $K_3$'s eigenvectors.

It is common sense that partial derivatives of a smooth field along any two of $x$-,$y$- and $z$-axis can be exchanged. Then it is expected that the discrete partial derivative operators $C_1, C_2, C_3$ discussed in previous section should commute with each other. This is indeed true in the case of FCC lattice[20]. However, use of this fact was far from being enough in our opinion.

In this section, we will prove some key results which are equivalent to $C_\ell C_{\ell'} = C_{\ell'} C_\ell$, $\ell, \ell' = 1,2,3$, for the triclinic lattice. Particularly, this commutativity and structure of eigenvectors of a (block) companion matrix will play a central role in deriving important eigen-decompositions of $C_\ell$. With these apparatuses, the whole process of derivation turns out very elegant and reader-friendly.

**Lemma 1.** *Let $p(t) = c_0 + c_1 t + \cdots + c_{n-1}t^{n-1} + t^n$ be an $n$-th degree complex monic polynomial with its companion matrix*

$$C_F(p) = \begin{bmatrix} 0 & 1 & 0 & \cdots & 0 \\ 0 & 0 & 1 & \cdots & 0 \\ \vdots & \vdots & \vdots & \ddots & \vdots \\ 0 & 0 & 0 & \ddots & 1 \\ -c_0 & -c_1 & -c_2 & \cdots & -c_{n-1} \end{bmatrix}, \tag{38}$$

*then $p(\lambda) = \det(\lambda - C_F(p))$, and especially the eigenvector of $C_F(p)$ corresponding to eigenvalue $\lambda_j$ is $[1, \lambda_j, \lambda_j^2, \cdots, \lambda_j^{n-1}]^\top, j = 1,2,\cdots,n$.*

Since lemma 1 can be directly verified, we skip its proof. Letting $c_1 = \cdots = c_{n-1} = 0$ and $c_0 = -\xi(\mathbf{k}\cdot\mathbf{a}_1)$ in lemma 1, we have the following theorem.

**Theorem 1** ([20]). *The eigenpairs of $K_1$ in (16) are $(\xi(\theta_{\mathbf{a}_1})\xi(i/n_1),\ X_i)$, where*

$$X_i = \left[1, \xi(\theta_{\mathbf{a}_1})\xi\left(\frac{i}{n_1}\right), \cdots, \xi((n_1-1)\theta_{\mathbf{a}_1})\xi\left(\frac{(n_1-1)i}{n_1}\right)\right]^\top, \tag{39}$$

*and $\theta_{\mathbf{a}_1} = \mathbf{k}\cdot\mathbf{a}_1/n_1$, for $i = 1, \ldots, n_1$.*

**Lemma 2** ([13]). *Let $M(\lambda) = M_0 + M_1\lambda + \cdots + M_{n-1}\lambda^{n-1} + \lambda^n$ with $M_j \in \mathbb{C}^{n\times n}, j = 0, 1, \cdots, n-1$, then $\det M(\lambda) = \det(\lambda - C_{BF}(M))$ where*

$$C_{BF}(M) = \begin{bmatrix} 0 & I_m & 0 & \cdots & 0 \\ 0 & 0 & I_m & \cdots & 0 \\ \vdots & \vdots & \vdots & \ddots & \vdots \\ 0 & 0 & 0 & \ddots & I_m \\ -M_0 & -M_1 & -M_2 & \cdots & -M_{n-1} \end{bmatrix}. \tag{40}$$

*Particularly if $v \in \mathbb{C}^m$ and $\lambda_0 \in \mathbb{C}$ satisfy $M(\lambda_0)v = 0$, then the eigenvector of $C_{BF}(M)$ corresponding to eigenvalue $\lambda_0$ is $[1, \lambda_0, \lambda_0^2, \cdots, \lambda_0^{n-1}]^\top \otimes v$.*

Similarly, letting $M_1 = M_2 = \cdots = M_{n-1} = 0$ and $M_0 = -\xi(\mathbf{k}\cdot\mathbf{a}_2)J_2$ in lemma 2, we see that eigenpairs of $K_2$ in (22) are made from those of $J_2$ in (20). If $\lambda_0$ is an eigenvalue of $J_2$ then the $n_2$-th root of $\xi(\mathbf{k}\cdot\mathbf{a}_2)\lambda_0$ is that of $K_2$.

**Lemma 3.** *Given $0 \neq \theta \in \mathbb{R}$ and $n \in \mathbb{N}$, for any $q \in Ind = \{1, 2, \cdots, n-1\}$,*

$$G_n(\theta, q) \colon = \begin{bmatrix} 0 & I_{n-q} \\ \xi(\theta)I_q & 0 \end{bmatrix}_{n\times n} = G_n(\theta, 1)^q. \tag{41}$$

*Proof.* When $q = 1$, (41) is obviously correct. Suppose (41) is correct when $1 \leq q = r < n-2$, *i.e.*, $G_n(\theta, r) = G_n(\theta, 1)^r$, then by direct multiplication we have

$$G_n(\theta, r)G_n(\theta, 1) = \begin{bmatrix} 0 & I_{n-r-1} \\ \xi(\theta)I_{r+1} & 0 \end{bmatrix} = G_n(\theta, r+1) = G_n(\theta, 1)^{r+1}. \tag{42}$$

By induction, we know that (41) is correct for all $q \in$ Ind. □

**Corollary 1.** *Given nonzero numbers $\theta \in \mathbb{R}$ and $n, q \in \mathbb{N}, n > q$, matrices*

$$G_n(\theta, 1) = \begin{bmatrix} 0 & I_{n-1} \\ \xi(\theta) & 0 \end{bmatrix},\ M = \begin{bmatrix} 0 & \xi(-\theta)I_q \\ I_{n-q} & 0 \end{bmatrix}$$

*commute. Let $\gamma_i = \xi(i/n)\xi(\theta/n), i = 1, 2, \ldots, n$, then eigenpairs of $M$ are $(\gamma_i^{-q}, v_i)$, where $v_i = [1, \gamma_i, \gamma_i^2, \cdots, \gamma_i^{n-1}]^\top$.*

*Proof.* Note that $M^* = G_n(\theta, q)$ and that eigenpairs of $G_n(\theta, 1)$ are $(\gamma_i, v_i)$. Hence by lemma 3 eigenpairs of $M^* = G_n(\theta, 1)^q$ are $(\gamma_i^q, v_i)$. Yet $MM^* = I_n$, therefore $Mv_i = v_i\gamma_i^{-q}$. Then it is obvious that $G_n(\theta, 1)M = MG_n(\theta, 1)$. □

Now let $\theta = \mathbf{k}\cdot\mathbf{a}_1,\ q = m_1,\ n = n_1$ in corollary 1, then $G_n(\theta, 1)$ becomes $K_1$ in (16) and $M$ becomes $J_2$ in (20). By simple manipulation we have the following result about $K_2$, according to lemma 2 and corollary 1.

**Theorem 2.** *The eigenpairs of $K_2$ in (22) are $(\xi(\theta_{\hat{\mathbf{a}}_2,i})\xi(j/n_2), Y_{ij}\otimes X_i)$, where $X_i$ is stated in (39) and*

$$\theta_{\hat{\mathbf{a}}_2,i} = \frac{1}{n_2}\left(\mathbf{k}\cdot\left(\mathbf{a}_2 - \frac{m_1}{n_1}\mathbf{a}_1\right) - \frac{im_1}{n_1}\right), \tag{43a}$$

$$Y_{ij} = \left[1, \xi(\theta_{\hat{\mathbf{a}}_2,i})\xi\left(\frac{j}{n_2}\right), \cdots, \xi((n_2-1)\theta_{\hat{\mathbf{a}}_2,i})\xi\left(\frac{(n_2-1)j}{n_2}\right)\right]^\top. \tag{43b}$$

*for $i = 1, \cdots, n_1,\ j = 1, \cdots, n_2$.*

In order to deal with $K_3$ in (31) in the same way, we need to establish commutativity between $J_3$ in (32) and $K_2$ in (22).

**Lemma 4.** *Given nonzero numbers $\theta_1 \neq \theta_2 \in \mathbb{R}$ and $n_1, n_2, r, q_1, q_2 \in \mathbb{N}, 1 \leq r < n_2, 1 \leq q_1 < q_2 < n_1$, the following two matrices*

$$W_1 = \begin{bmatrix} 0 & \xi(-\theta_2)I_r \otimes G^\diamond_{n_1}(\theta_1, q_1) \\ I_{n_2-r}\otimes G^\diamond_{n_1}(\theta_1, q_2) & 0 \end{bmatrix},$$
$$W_2 = \begin{bmatrix} 0 & I_{n_2-1}\otimes I_{n_1} \\ \xi(\theta_2)G^\diamond_{n_1}(\theta_1, q_2-q_1) & 0 \end{bmatrix}$$

*commute, where $G^\diamond$ refers to $G$ and $G^*$, respectively.*

*Proof.* To avoid ambiguity, it is required that once $G^\diamond_{n_1}(\theta_1, q_1) = G_{n_1}(\theta_1, q_1)$, then $G^\diamond$ must stands for $G$ in all three places in $W_1, W_2$. So does $G^*$. By lemma 3, $G^\diamond_{n_1}(\theta_1, q_2-q_1)G^\diamond_{n_1}(\theta_1, q_1) = G^\diamond_{n_1}(\theta_1, q_2) = G^\diamond_{n_1}(\theta_1, q_1)G^\diamond_{n_1}(\theta_1, q_2-q_1)$. By direct block matrix multiplication, we can easily have

$$W_1W_2 = \begin{bmatrix} 0 & \xi(-\theta_2)I_{r-1}\otimes G^\diamond_{n_1}(\theta_1, q_1) \\ I_{n_2-r+1}\otimes G^\diamond_{n_1}(\theta_1, q_2) & 0 \end{bmatrix} = W_2W_1.$$

□

**Theorem 3.** *The eigenpairs of $K_3$ in (31) are $(\xi(\theta_{\hat{\mathbf{a}}_3,ij})\xi(k/n_3), Z_{ijk}\otimes Y_{ij}\otimes X_i)$, with $X_i$ and $Y_{ij}$ stated in (39) and (43b), respectively, where*

$$\theta_{\hat{\mathbf{a}}_3,ij} = \frac{1}{n_3}\left\{\mathbf{k}\cdot\hat{\mathbf{a}}_3 - \frac{m_3}{n_2}j + \left(\frac{m_1}{n_1}\frac{m_3}{n_2} - \frac{m_2}{n_1}\right)i\right\}, \tag{44a}$$

$$Z_{ijk} = \left[1, \xi(\theta_{\hat{\mathbf{a}}_3,ij})\xi\left(\frac{k}{n_3}\right), \cdots, \xi((n_3-1)\theta_{\hat{\mathbf{a}}_3,ij})\xi\left(\frac{(n_3-1)k}{n_3}\right)\right]^\top, \tag{44b}$$

$$\hat{\mathbf{a}}_3 = \mathbf{a}_3 - \frac{m_3}{n_2}\mathbf{a}_2 + \left(\frac{m_1}{n_1}\frac{m_3}{n_2} - \frac{m_2}{n_1}\right)\mathbf{a}_1,$$

*for $i = 1, \cdots, n_1,\ j = 1, \cdots, n_2,\ k = 1, \cdots, n_3$.*

*Proof.* In lemma 4 let $\theta_1 = \mathbf{k}\cdot\mathbf{a}_1,\ \theta_2 = \mathbf{k}\cdot\mathbf{a}_2,\ r = m_3,\ q_1 = m_2 - m_1,\ q_2 = m_2,\ G^\diamond = G^*$, then $J_3$ commutes with $K_2$. This implies that $J_3$ also has eigenvectors $v_{ij} = Y_{ij}\otimes X_i$ stated in theorem 2. If $\mu_{ij}$ is the eigenvalue of $J_3$ corresponding to $v_{ij}$, then comparing $(m_3n_1 + m_2 + 1)$-th entry of $J_3v_{ij}$ and $\mu_{ij}v_{ij}$ we have $\mu_{ij}\xi\left(m_3\theta_{\hat{\mathbf{a}}_2,i}\right)\xi(jm_3/n_2)\xi(m_2\theta_{\mathbf{a}_1})\xi(im_2/n_1) = 1$. Therefore corresponding to $v_{ij}$ the eigenvalue of $\xi(\mathbf{k}\cdot\mathbf{a}_3)J_3$ is

$$\xi(\mathbf{k}\cdot\mathbf{a}_3)\xi\left(-m_3\theta_{\hat{\mathbf{a}}_2,i}\right)\xi\left(-\frac{jm_3}{n_2}\right)\xi\left(-m_2\theta_{\mathbf{a}_1}\right)\xi\left(-\frac{im_2}{n_1}\right) = \xi(n_3\theta_{\hat{\mathbf{a}}_3,ij}).$$

Then, the $n_3$-th root of $\xi(n_3\theta_{\hat{\mathbf{a}}_3,ij})$ is just $\xi(\theta_{\hat{\mathbf{a}}_3,ij})\xi(k/n_3), k = 1,\cdots,n_3$. □

As before, due to different combination of $\phi_3,\ \phi_2,\ \ell_2$ different expressions of $\theta_{\hat{\mathbf{a}}_2,i},\ \hat{\mathbf{a}}_3,\ \theta_{\hat{\mathbf{a}}_3,ij}$ can be found in the Appendix.

Now we summarize the results obtained in this section.

We put all eigenvectors of $K_i$ in a matrix $T$,

$$T = [T_1,\quad T_2,\quad \cdots,\quad T_{n_1}]/\sqrt{n} \in \mathbb{C}^{n\times n}, \tag{45a}$$

$$T_i = [T_{i1},\quad T_{i2},\quad \cdots,\quad T_{in_2}] \in \mathbb{C}^{n\times(n_2n_3)}, \tag{45b}$$

$$T_{ij} = [Z_{ij1}\otimes Y_{ij}\otimes X_i,\ Z_{ij2}\otimes Y_{ij}\otimes X_i,\ \cdots,\ Z_{ijn_3}\otimes Y_{ij}\otimes X_i], \tag{45c}$$

for $i = 1,\cdots,n_1, j = 1,\cdots,n_2$. This eigenmatrix $T$ is unitary since it is straightforward to show that $K_1, K_2, K_3$ are normal (in fact unitary) matrices and each vector of $T$ is normalized.

For the eigenvalues, we set

$$\begin{aligned}
\Lambda_{n_1} &= \operatorname{diag}\left(\xi(\theta_{\mathbf{a}_1})\xi([1:n_1]^\top/n_1) - 1\right)/\delta_x, \quad \Lambda_1 = \Lambda_{n_1}\otimes I_{n_2}\otimes I_{n_3}\\
\Lambda_{in_2} &= \operatorname{diag}\left(\xi(\theta_{\hat{\mathbf{a}}_2,i})\xi([1:n_2]^\top/n_2) - 1\right)/\delta_y, \quad \Lambda_2 = \oplus_{i=1}^{n_1}(\Lambda_{in_2}\otimes I_{n_3}),\\
\Lambda_{ijn_3} &= \operatorname{diag}\left(\xi(\theta_{\hat{\mathbf{a}}_3,ij})\xi([1:n_3]^\top/n_3) - 1\right)/\delta_z, \quad \Lambda_3 = \oplus_{i=1}^{n_1}\oplus_{j=1}^{n_2}\Lambda_{ijn_3}.
\end{aligned}$$

Then, from Theorems 1, 2, and 3, it holds that

$$C_1T = T\Lambda_1,\quad C_2T = T\Lambda_2,\quad C_3T = T\Lambda_3. \tag{46}$$

Recall that in [11, 20], we have derived the eigen-decompositions (46) only for SC and FCC lattices. Now it is clear that the formalism is the same for all Bravais lattices, though $\theta_{\mathbf{a}_1}$, $\theta_{\hat{\mathbf{a}}_2,i}$ and $\theta_{\hat{\mathbf{a}}_3,ij}$ depend on $\mathbf{a}_1$, $\mathbf{a}_2$, $\mathbf{a}_3$.

## 6. Eigen-decomposition of $\mathcal{A}$

On basis of the results of previous section, we can proceed to derive eigen-decomposition of $\mathcal{A}$, especially to obtain the range-space of $\mathcal{A}$ explicitly. Rather than find singular value decomposition (SVD) of $\mathcal{C}$ as is done in Ref. [20, 21], we will preserve the complex skew-symmetry which is intrinsic to $\mathcal{C}$.

From Eq. (46), we know that $\mathcal{C}$ in (37) is unitarily similar to a complex skew-symmetric matrix $\mathbf{\Lambda}$

$$\mathbf{\Lambda} = \begin{bmatrix} 0 & -\Lambda_3 & \Lambda_2 \\ \Lambda_3 & 0 & -\Lambda_1 \\ -\Lambda_2 & \Lambda_1 & 0 \end{bmatrix} = (I_3 \otimes T)^* \mathcal{C} (I_3 \otimes T) \tag{47}$$

Moreover, doing a perfect shuffle of this $\Lambda$, *i.e.*, multiplying

$$P = [e_1, e_{n+1}, e_{2n+1}, e_2, e_{n+2}, e_{2n+2}, \cdots, e_n, e_{2n}, e_{3n}] \in \mathbb{R}^{3n \times 3n}, \tag{48}$$

from the right side and $P^\top$ from the left side, we can transform $\Lambda$ to a block diagonal matrix

$$P^\top \Lambda P = L_1 \oplus L_2 \oplus \cdots \oplus L_n,$$

with $L_\ell \in \mathbb{C}^{3\times 3}, L_\ell^\top = -L_\ell, \ell = 1, 2, \cdots, n$. That means we can just deal with each block $L_\ell$ separately. SVD in this case does not reveal the underlying structure, therefore is not preferred here. Yet it is well-known that [19] canonical form of a complex skew-symmetric matrix under unitary congruence is a real quasi-diagonal skew-symmetric matrix, which is certainly rank-revealing. In particular, we find that for the $3 \times 3$ complex skew-symmetric matrix, the canonical form can be done in almost one step, as shown below, which is simpler than some well-established algorithm to achieve the same goal. Then we can express analytically the range space of a $L_\ell$, which is of rank 2.

The following lemma is a generalization of the scalar triple product in $\mathbb{R}^3$.

**Lemma 5.** *Given a nonzero vector $c = [c_1, c_2, c_3]^\top \in \mathbb{C}^3$, we define a corresponding skew-symmetric matrix*

$$L = \begin{bmatrix} 0 & -c_3 & c_2 \\ c_3 & 0 & -c_1 \\ -c_2 & c_1 & 0 \end{bmatrix}_{3\times 3}, L^\top = -L.$$

*Then $Lc = [0,0,0]^\top$ and for any $u = [u_1, u_2, u_3]^\top, v = [v_1, v_2, v_3]^\top \in \mathbb{C}^3$,*

$$v^\top L u = \det \left( \begin{bmatrix} v_1 & v_2 & v_3 \\ c_1 & c_2 & c_3 \\ u_1 & u_2 & u_3 \end{bmatrix} \right) = \det([v, c, u]).$$

The next lemma follows from one basic algorithm in numerical linear algebra.

**Lemma 6.** *Given a nonzero vector $c = [c_1, c_2, c_3]^\top \in \mathbb{C}^3$, the following Householder matrix*

$$H = I_3 - \tau [1, h_2, h_3]^\top \left[1, \overline{h_2}, \overline{h_3}\right]$$

*satisfies*

$$H^* H = I_3, H^* c = [\beta, 0, 0]^\top,$$

*where*

$$\beta = -sign(\Re(c_1))||c||_2,\ \tau = \frac{\beta - c_1}{\beta},\ h_\ell = \frac{c_\ell}{c_1 - \beta},\ \ell = 2, 3. \tag{49}$$

*Furthermore,* $\det(H) = (\beta - c_1)/(\overline{c_1} - \beta)$.

*Proof.* Here we only prove the last equality, since the rest is just the basic algorithm in LAPACK[12] to compute the complex Householder matrix. By Sylvester's determinant identity[5], we have

$$\det(H) = 1 - \tau||[1, h_2, h_3]||_2 = 1 - \frac{\beta - c_1}{\beta} - \frac{\beta - c_1}{\beta}\frac{|c_2|^2 + |c_3|^2}{|c_1 - \beta|^2}$$
$$= \frac{c_1}{\beta} + \frac{|c_2|^2 + |c_3|^2}{\beta(\overline{c_1} - \beta)} = \frac{|c_1|^2 + |c_2|^2 + |c_3|^2 - \beta c_1}{\beta(\overline{c_1} - \beta)} = \frac{\beta - c_1}{\overline{c_1} - \beta}.$$

□

**Theorem 4.** *Given a nonzero vector $c = [c_1, c_2, c_3]^\top \in \mathbb{C}^3$, $L$ defined in lemma 5 and $H$ defined in lemma 6, let $V = H\,diag(1, 1, \det(H^*))$, then*

$$V^\top LV = \begin{bmatrix} 0 & 0 & 0 \\ 0 & 0 & -\beta \\ 0 & \beta & 0 \end{bmatrix}, \quad V^*V = VV^* = I_3. \tag{50}$$

*That is to say,*

$$L = \beta\overline{V}(:, 3)\overline{V}(:, 2)^\top - \beta\overline{V}(:, 2)\overline{V}(:, 3)^\top, \tag{51a}$$

$$\overline{V}(:, 2) = \overline{V(:, 2)} = \left[\frac{c_2(\overline{c_1} - \beta)}{\beta(c_1 - \beta)}, 1 + \frac{|c_2|^2}{\beta(c_1 - \beta)}, \frac{\overline{c_3}c_2}{\beta(c_1 - \beta)}\right]^\top, \tag{51b}$$

$$\overline{V}(:, 3) = \overline{V(:, 3)} = \left[-\frac{c_3}{\beta}, \frac{\overline{c_2}c_3}{\beta(\beta - \overline{c_1})}, \frac{c_1}{\beta} - \frac{|c_2|^2}{\beta(\beta - \overline{c_1})}\right]^\top. \tag{51c}$$

*Proof.* Let us consider the skew-symmetric matrix $H^\top LH$ first. Obviously diagonal entries of $H^\top LH$ vanish. Since $H^*c = \beta e_1$, then $He_1 = c/\beta$, therefore $LHe_1 = [0, 0, 0]^\top$ and $e_1^\top H^\top L = [0, 0, 0]$. That means $e_3^\top(H^\top LH)e_2 = -e_2^\top(H^\top LH)e_3 \neq 0$. In fact, we have

$$e_3^\top(H^\top LH)e_2 = (He_3)^\top L(He_2)$$
$$= \det([He_3, c, He_2]) = \det(H)\det([e_3, H^*c, e_2])$$
$$= \det(H)\det([e_3, \beta e_1, e_2]) = \beta\det(H).$$

Since $\det(H^*)\det(H) = 1$, we have $e_3^\top(V^\top LV)e_2 = \beta = -e_2^\top(V^\top LV)e_3$.

Here we calculate $\overline{V}(:, 3)$ explicitly, which is somewhat tedious.

$$\overline{V}(1, 3) = -\det(H)\overline{\tau}h_3 = \frac{\beta - c_1}{\overline{c_1} - \beta}\frac{c_3(\overline{c_1} - \beta)}{\beta(c_1 - \beta)} = -\frac{c_3}{\beta},$$

$$\overline{V}(2, 3) = -\det(H)\overline{\tau}h_3\overline{h_2} = \frac{\beta - c_1}{\overline{c_1} - \beta}\frac{\overline{c_1} - \beta}{\beta}\frac{\overline{c_2}c_3}{|c_1 - \beta|^2} = \frac{\overline{c_2}c_3}{\beta(\beta - \overline{c_1})},$$

$$\overline{V}(3, 3) = \det(H) - \det(H)\overline{\tau}|h_3|^2 = \frac{\beta - c_1}{\overline{c_1} - \beta} + \frac{\beta - c_1}{\overline{c_1} - \beta}\frac{\overline{c_1} - \beta}{\beta}\frac{|c_3|^2}{|c_1 - \beta|^2}$$
$$= \frac{\beta(\beta - c_1) - |c_3|^2}{\beta(\overline{c_1} - \beta)} = \frac{|c_1|^2 + |c_2|^2 - \beta c_1}{\beta(\overline{c_1} - \beta)} = \frac{c_1}{\beta} + \frac{|c_2|^2}{\beta(\overline{c_1} - \beta)}.$$

$\overline{V}(:, 2)$ can be similarly calculated and is skipped. □

In fact, letting $\widetilde{V} = V(:,[2,3]) \in \mathbb{C}^{3\times 2}$, the canonical form of $L$ in theorem 4 can be reduced to the following, without keeping the nullspace

$$L = \overline{\widetilde{V}} \begin{bmatrix} 0 & -\beta \\ \beta & 0 \end{bmatrix} \widetilde{V}^*, \quad \widetilde{V}^*\widetilde{V} = I_2 = \widetilde{V}^\top \overline{\widetilde{V}},$$

and similarly for each $L_\ell, \ell = 1, 2, \cdots, n$, we have

$$L_\ell = \overline{\widetilde{V}}_\ell \begin{bmatrix} 0 & -\beta_\ell \\ \beta_\ell & 0 \end{bmatrix} \widetilde{V}_\ell^*, \quad \widetilde{V}_\ell^*\widetilde{V}_\ell = I_2 = \widetilde{V}_\ell^\top \overline{\widetilde{V}}_\ell.$$

Finally the eigenspace of $\mathcal{A}$ which is orthogonal to the nullspace of $\mathcal{A}$ can be derived using key results in this and previous section[20, Theorem 3.7].

**Theorem 5.** *Let $\mathcal{A}$, $T$, $P$ defined in* (36),(45a),(48) *respectively, and denote*

$$\begin{aligned}
\mathcal{V}_r &:= blkdiag\left(\widetilde{V}_1, \widetilde{V}_2, \cdots, \widetilde{V}_n\right) \in \mathbb{C}^{3n\times 2n}, \\
\Lambda_r &:= diag\left(\beta_1^2, \beta_1^2, \beta_2^2, \beta_2^2, \cdots, \beta_n^2, \beta_n^2\right) \in \mathbb{C}^{2n\times 2n}, \\
\mathcal{Q}_r &:= \left(I_3 \otimes T\right) P \mathcal{V}_r \in \mathbb{C}^{3n\times 2n},
\end{aligned}$$

*then*

$$\mathcal{A} = \mathcal{Q}_r \Lambda_r \mathcal{Q}_r^*, \quad \mathcal{Q}_r^* \mathcal{Q}_r = I_{2n}. \tag{52}$$

## 7. Iterative solver FAME for NFSEP (53)

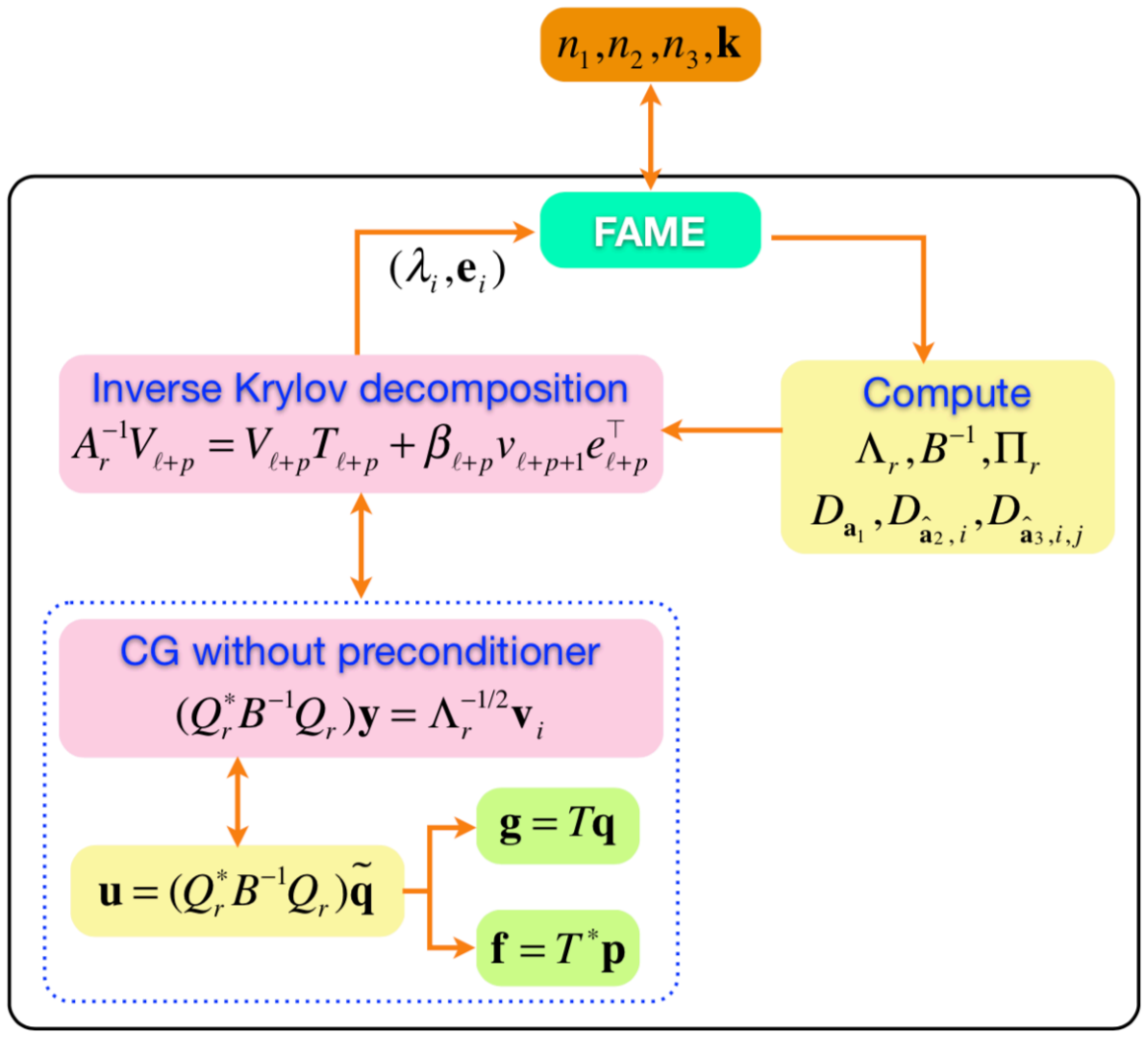

Figure 4: Flow charts of FAME for NFSEP (53).

Eventually, all previous derivations show that the nullspace free method proposed in [20] still works for all Bravais lattice, which transforms the GEP (36) into the following NFSEP:

$$\mathcal{A}_r\widetilde{E} = \lambda\widetilde{E}, \quad E = \mathcal{B}^{-1}\mathcal{Q}_r\Lambda_r^{1/2}\widetilde{E}, \quad \mathcal{A}_r = \Lambda_r^{1/2}\mathcal{Q}_r^*\mathcal{B}^{-1}\mathcal{Q}_r\Lambda_r^{1/2} = \mathcal{A}_r^* > 0. \quad (53)$$

Since the nullspace of the GEP (36) has been completely deflated, the first challenge mentioned in Sec. 1 is resolved.

To solve (53), a fast eigensolver called FAME was proposed in [20] originally for SC and FCC lattices, and can also be similarly applied to all Bravais lattices. The flowchart of our fast eigensolver FAME is shown in Figure 4. As shown in this figure, conjugate gradient (CG) method without preconditioner to solve the linear system is very efficient, because the condition number of $\mathcal{Q}_r^*\mathcal{B}^{-1}\mathcal{Q}_r$ in (53) is bounded by that of $\mathcal{B}^{-1}$. The second challenge mentioned in Sec. 1 is resolved in the case of positive diagonal $\mathcal{B}$.

In CG method, multiplying a column vector $\widetilde{\mathbf{q}}$ by $\mathcal{Q}_r^*\mathcal{B}^{-1}\mathcal{Q}_r$ is essentially reduced to $T\mathbf{q}$ and $T^*\mathbf{p}$ besides some diagonal scalings, where $\mathbf{q}, \mathbf{p}$ are intermediate variables. Fortunately, we discover that the most expensive operations $T\mathbf{q}$ and $T^*\mathbf{p}$ can be efficiently computed via Algorithms 1 and 2, respectively, which are described in Appendix E. In a nutshell, these two algorithms are just wrappers for backward and forward FFT, respectively, harnessing (30).

## 8. Numerical Experiments

To demonstrate the accuracy and efficiency of our framework, we call functions **eigs**, **pcg**, **fft** and **ifft** of Matlab[3] R2017b to implement key operations in our fast eigensolver FAME and calculate band strucuture of one benchmark system of the double gyroid PhC[28] in Body-Centered Cubic (BCC) lattice. In our calculation, the convergence tolerance of **eigs** and **pcg** is set to $10^{-12}$ and $10^{-13}$, respectively. All computations are performed on an Intel (R) Xeon (R) E5-2643 3.30GHz processor with 96 GB RAM in double precision arithmetic.

The lattice translation vectors $\mathbf{a}_1, \mathbf{a}_2, \mathbf{a}_3$ of the BCC lattice are

$$\mathbf{a}_1 = \frac{a}{2}\left[-1,1,1\right]^\top, \ \mathbf{a}_2 = \frac{a}{2}\left[1,-1,1\right]^\top, \ \mathbf{a}_3 = \frac{a}{2}\left[1,1,-1\right]^\top,$$

where $a$ is the lattice constant. The reciprocal lattice vectors $\mathbf{b}_1, \mathbf{b}_2, \mathbf{b}_3$ satisfy $[\mathbf{b}_1\ \mathbf{b}_2\ \mathbf{b}_3]\left[\mathbf{a}_1\ \mathbf{a}_2\ \mathbf{a}_3\right]^\top = 2\pi I_3$. The vertexes of the Brillouin zone (see Figure 5(b)) of BCC lattice can be represented in the oblique coordinate system spanned by $\mathbf{b}_1$, $\mathbf{b}_2$ and $\mathbf{b}_3$ as

$$\Gamma = [0,0,0]^\top, \ H = \left[\frac{1}{2}, -\frac{1}{2}, \frac{1}{2}\right]^\top, P = \left[\frac{1}{4}, \frac{1}{4}, \frac{1}{4}\right]^\top, N = \left[0, \frac{1}{2}, 0\right]^\top, H' = \left[-\frac{1}{2}, \frac{1}{2}, \frac{1}{2}\right]^\top.$$

Let $\mathbf{r} = (x, y, z)$. The double gyroid region in Figure 5(a) can be described by the set $\{r \in \mathbb{R}^3 | g(\mathbf{r}) > 1.1\} \cup \{r \in \mathbb{R}^3 | g(-\mathbf{r}) > 1.1\}$ where $g(\mathbf{r}) = \sin(2\pi x/a)\cos(2\pi y/a) + \sin(2\pi y/a)\cos(2\pi z/a) + \sin(2\pi z/a)\cos(2\pi x/a)$. For

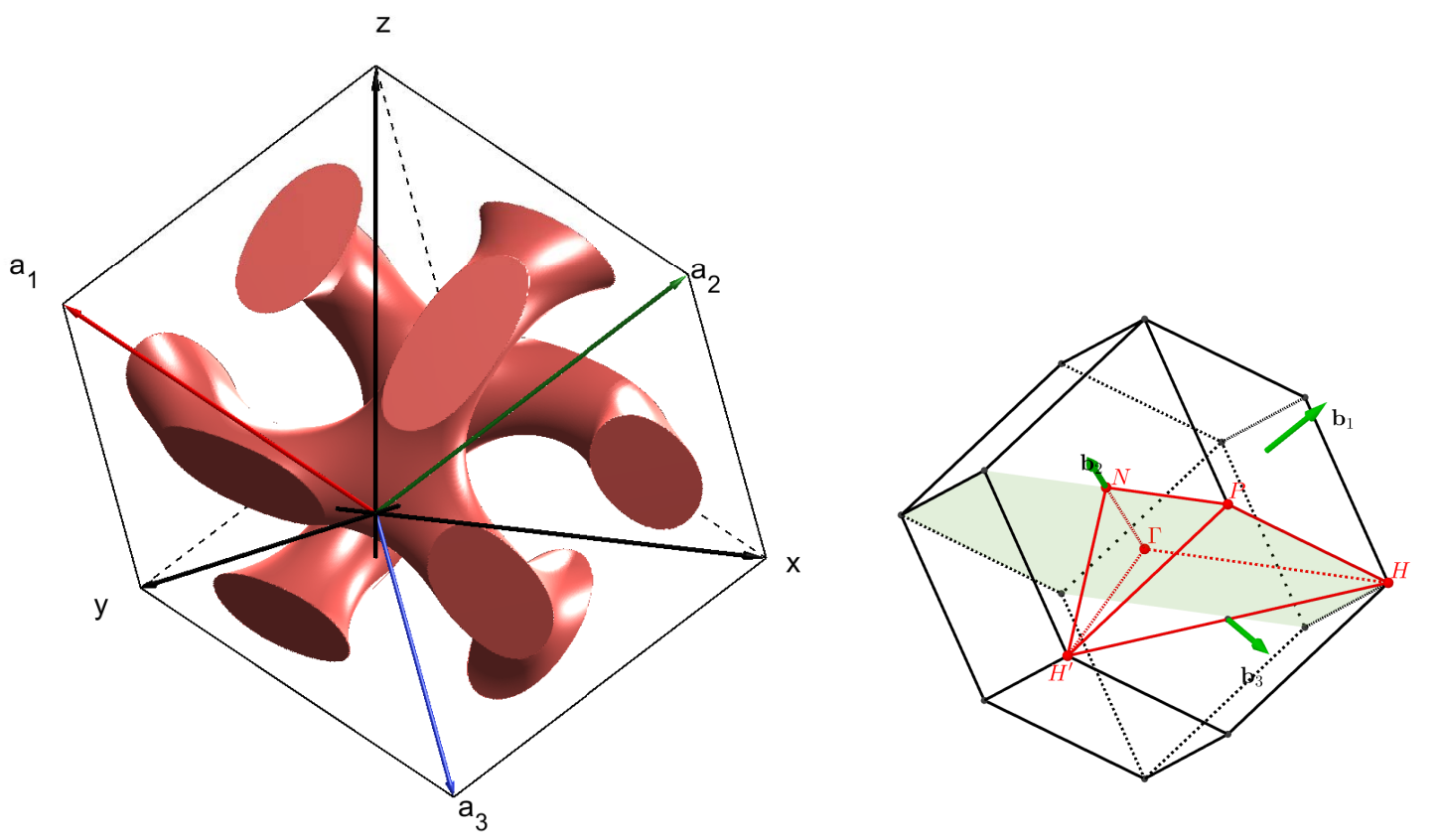

(a) The physical cell of BC lattice and double gyroid structure.

(b) Brillouin zone and its vertex of BCC lattice.

Figure 5: Illustration of one PhC in BC lattice and its Brillouin zone

convenience we set $a = 1$. Suppose that the double gyroid region is filled up material whose permittivity is uniformly $\varepsilon = 16$ and that the rest part is just vacuum ($\varepsilon = 1$), then we compute ten smallest positive eigenvalues and associated eigenvectors of the this system.

The band structure shown in Figure 6(a) does not show any discernible discrepancy with the one in Ref. [28], which partially evidences the accuracy of our method. Even the dimension of the NFSEP (53) is as large as $3,456,000$, it takes at most $7 \times 10^3$ seconds to calculate ten target eigenpairs at each **k**-point as shown in Figure 6(b) (1), which is acceptable considering the serials implementation. More detailedly, in Figure 6(b) (2) the number of iterations in function **eigs** versus **k** is plotted, where we can see that the Invert-Lanczos process converges in 60 to 170 steps for the ten target eigenpairs given **k**. In Figure 6(b) (3), the number of iterations in function **pcg** without preconditioner versus **k** is plotted, where we can see that on average we need 34 to 42 iterations to solve the linear system in one step of Invert-Lanczos process. The overall efficiency of our algorithm is impressive.

## 9. Conclusion

In a word, the major contribution we have made in the present work is the establishment of a complete and unified framework to solve Maxwells Eigenvalue Problem for 3D isotropic photonic crystals in all 14 Bravais lattices. IT is highlighted that our method is remarkably efficient. Compared with $\mathcal{O}(n^2)$ of other method, the overall computational complexity of our method is $\mathcal{O}(n \log n)$, thanks to the feasibility of FFT algorithm in our framework, which is actually rooted in the eigen-decomposition of discrete partial derivative operators

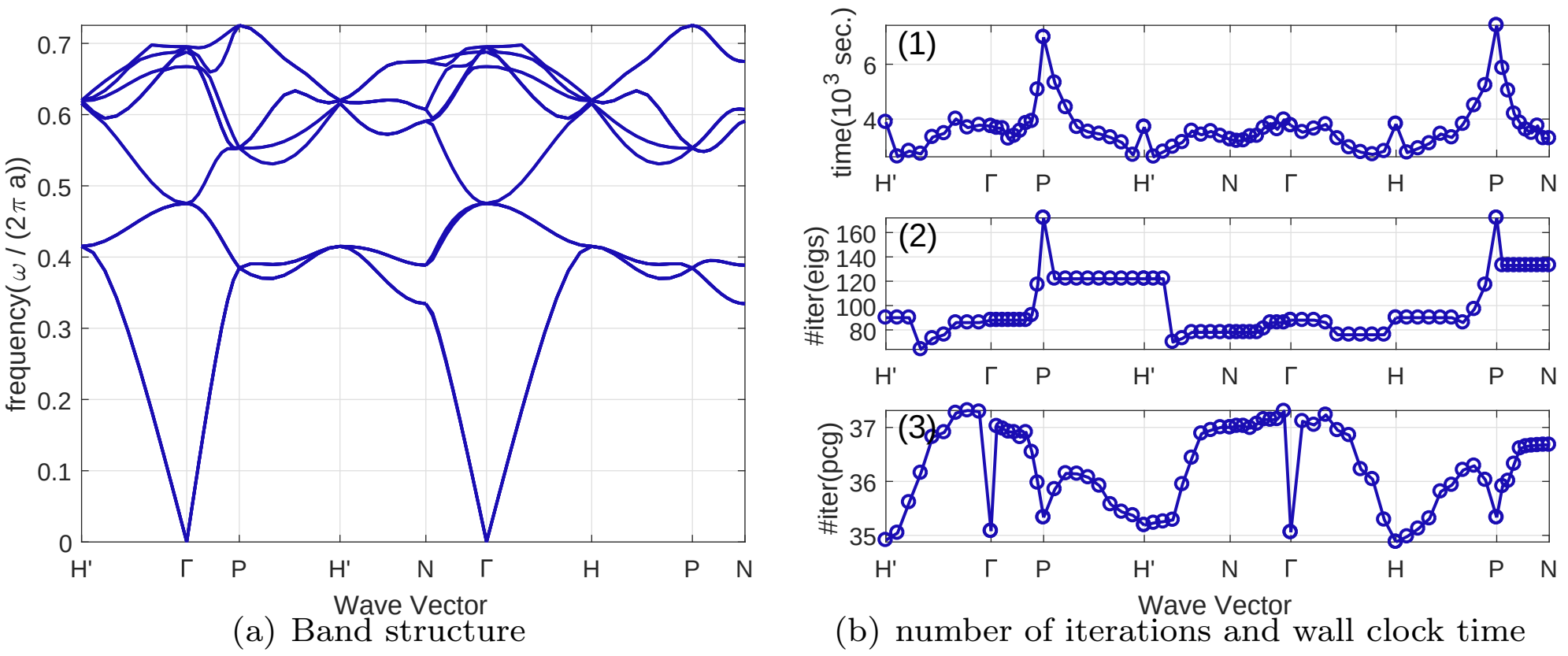

(a) Band structure

(b) number of iterations and wall clock time

Figure 6: (a) the band structure of double gyroid PhC; (b)(1) the average number of iterations in **pcg** without preconditioner, (b)(2) the number of iterations in **eigs**, (b)(3) wall clock time spent on ten target eigenvalues.

$\partial_x, \partial_y, \partial_z$ with reformulated Bloch-boundary condition. The commutativity among discrete partial derivative operators is one of the key machinery that allows us to derive these eigen-decompositions in a light way. On the other hand, the fast convergence of our eigensolver FAME is guaranteed by the novel nullspace free method that thoroughly removes the considerable nullspace of the discrete double-curl operator $\mathcal{A}$.

Also, our unique way to compute the canonical form of a $3 \times 3$ complex skew-symmetric matrix under unitary congruence may be of independent interest. The significance of the 3D cubic working cell defined in Sec. 2 of this work will be discussed exhaustively elsewhere. Extension of our present framework to 3D anisotropic photonic crystals is under investigation and will be reported in near future.

## Acknowledgments

T.-M. Huang was partially supported by the Ministry of Science and Technology (MoST) 105-2115-M-003-009-MY3, National Centre of Theoretical Sciences (NCTS) in Taiwan. T. Li was supported in parts by the NSFC 11471074. W.-W. Lin was partially supported by MoST 106-2628-M-009-004-, NCTS and ST Yau Centre in Taiwan. Prof. So-Hsiang Chou is greatly appreciated for his valuable feedback and suggestion on this manuscript.

## Appendix A. Coordinates of lattice vectors of 7 lattice systems

In Table. A.1, we list the Cartesian coordinates of lattice vectors $\tilde{\mathbf{a}}_1, \tilde{\mathbf{a}}_2, \tilde{\mathbf{a}}_3$ of all 7 lattice systems.

| | Triclinic | Monoclinic |
|---|---|---|
| Primitive (Pr.) | $\begin{bmatrix} a_1 & a_2\cos\phi_3 & a_3\cos\phi_2 \\ 0 & a_2\sin\phi_3 & a_3\ell_2 \\ 0 & 0 & a_3\ell_3 \end{bmatrix}$ | $\begin{bmatrix} a_1 & a_2\cos\phi_3 & 0 \\ 0 & a_2\sin\phi_3 & 0 \\ 0 & 0 & a_3 \end{bmatrix}$ |
| A-Base-centered (A-Ba.C.) | | $\frac{1}{2}\begin{bmatrix} 2a_1 & a_2\cos\phi_3 & a_2\cos\phi_3 \\ 0 & a_2\sin\phi_3 & a_2\sin\phi_3 \\ 0 & a_3 & -a_3 \end{bmatrix}$, if $a_1 \geq \frac{\sqrt{a_2^2+a_3^2}}{2}$<br>$\frac{1}{2}\begin{bmatrix} a_2\cos\phi_3 & a_2\cos\phi_3 & 2a_1 \\ a_2\sin\phi_3 & a_2\sin\phi_3 & 0 \\ a_3 & -a_3 & 0 \end{bmatrix}$, otherwise |

| | Rhombohedral | Cubic |
|---|---|---|
| Primitive | | $\begin{bmatrix} a_1 & 0 & 0 \\ 0 & a_1 & 0 \\ 0 & 0 & a_1 \end{bmatrix}$ |
| Body-centered (B.C.) | | $\frac{a_1}{2}\begin{bmatrix} -1 & 1 & 1 \\ 1 & -1 & 1 \\ 1 & 1 & -1 \end{bmatrix}$ |
| Face-centered (F.C.) | | $\frac{a_1}{2}\begin{bmatrix} 1 & 0 & 1 \\ 1 & 1 & 0 \\ 0 & 1 & 1 \end{bmatrix}$ |
| Rhombohedrally-centered (R.C.) | $\begin{bmatrix} 0 & a_1/2 & -a_1/2 \\ -a_1/\sqrt{3} & \sqrt{3}a_1/6 & \sqrt{3}a_1/6 \\ a_3/3 & a_3/3 & a_3/3 \end{bmatrix}$ | |

| | Orthorhombic | Tetragonal | Hexagonal |
|---|---|---|---|
| Primitive | $\begin{bmatrix} a_1 & 0 & 0 \\ 0 & a_2 & 0 \\ 0 & 0 & a_3 \end{bmatrix}$ | $\begin{bmatrix} a_1 & 0 & 0 \\ 0 & a_1 & 0 \\ 0 & 0 & a_3 \end{bmatrix}$ | $\begin{bmatrix} a_1 & -a_1/2 & 0 \\ 0 & \sqrt{3}a_1/2 & 0 \\ 0 & 0 & a_3 \end{bmatrix}$ |
| A-Base-centered | $\frac{1}{2}\begin{bmatrix} 2a_1 & 0 & 0 \\ 0 & a_2 & a_2 \\ 0 & a_3 & -a_3 \end{bmatrix}$ | | |
| C-Base-centered (C-Ba.C.) | $\frac{1}{2}\begin{bmatrix} a_1 & -a_1 & 0 \\ a_2 & a_2 & 0 \\ 0 & 0 & 2a_3 \end{bmatrix}$ | | |
| Body-centered | $\frac{1}{2}\begin{bmatrix} -a_1 & a_1 & a_1 \\ a_2 & -a_2 & a_2 \\ a_3 & a_3 & -a_3 \end{bmatrix}$ | $\frac{1}{2}\begin{bmatrix} a_1 & -a_1 & a_1 \\ a_1 & a_1 & -a_1 \\ -a_3 & a_3 & a_3 \end{bmatrix}$ | |
| Face-centered | $\frac{1}{2}\begin{bmatrix} a_1 & a_1 & 0 \\ a_2 & 0 & a_2 \\ 0 & a_3 & a_3 \end{bmatrix}$ | | |

Table A.1: Lattice vectors $[\tilde{\mathbf{a}}_1, \tilde{\mathbf{a}}_2, \tilde{\mathbf{a}}_3]$ for 7 lattice systems, with notations specified in Sec. 2.

## Appendix B. derivation Figure 2 and BC (11)

It is best to visualize the investigation starting from Figure 7(a), where we have $\phi_3 < \pi/2,\ \phi_2 < \pi/2,\ \ell_2 > 0$. Results of other possibilities such as $\phi_3 < \pi/2,\ \phi_2 > \pi/2,\ \ell_2 > 0$ will be discussed later.

In Figure 7(a), suppose $\square \mathrm{OR}_1\mathrm{R}_2\mathrm{R}_3$ is the bottom surface of $\mathbb{D}$, while $\square \mathrm{R}_4\mathrm{R}_5\mathrm{R}_6\mathrm{R}_7$ is image of the top surface of $\mathbb{D}$ under $\mathcal{T}_{-\mathbf{a}_3}$, which contains the origin in this case. Also, naturally we have the 2D oblique coordinate system with $\mathbf{a}_1$-,$\mathbf{a}_2$-axis. With slight abuse of notation, I,II,III,IV denote four patches of the $\square \mathrm{R}_4\mathrm{R}_5\mathrm{R}_6\mathrm{R}_7$, located in the first, second, third, fourth quadrant, respectively, of this oblique coordinate system. Our goal is to map $\square \mathrm{R}_4\mathrm{R}_5\mathrm{R}_6\mathrm{R}_7$ to $\square \mathrm{OR}_1\mathrm{R}_2\mathrm{R}_3$, respecting the periodicity along $\mathbf{a}_1, \mathbf{a}_2$. Here we have the 2D physical cell generated by $\mathbf{a}_1$ and $\mathbf{a}_2$, *i.e.*, the set $\{\alpha\mathbf{a}_1 + \beta\mathbf{a}_2 : \alpha, \beta \in [0, 1)\}$, and its periodic images under $\mathcal{T}_{\mathbf{a}_1}, \mathcal{T}_{\mathbf{a}_2}$ which fill up the whole plane, *i.e.*, the set $\{\alpha\mathbf{a}_1 + \beta\mathbf{a}_2 : \alpha, \beta \in \mathbb{R}\}$. Due to the periodicity, it is best to reduce all objects on the plane to their counterparts within the 2D physical cell. The rule is that whenever a point is outside the 2D physical cell, *i.e.*, $\alpha, \beta \notin [0, 1)$, we evaluate its image within the 2D physical cell under modulo operation. For example, for points in patch III we have $\alpha, \beta \in [-1, 0)$, then due to

$$\alpha\mathbf{a}_1 + \beta\mathbf{a}_2 \equiv (1 + \alpha)\mathbf{a}_1 + (1 + \beta)\mathbf{a}_2 = \mathcal{T}_{\mathbf{a}_1}\mathcal{T}_{\mathbf{a}_2}(\alpha\mathbf{a}_1 + \beta\mathbf{a}_2) \mod \mathbf{a}_1, \mathbf{a}_2,$$

patch III is mapped to its counterpart in the 2D physical cell shown in Figure 7(b). Other patches are similarly relocated.

As shown in Figure 7(c), it is easy to map the 2D physical cell to $\square \mathrm{OR}_1\mathrm{R}_2\mathrm{R}_3$, which is realized if triangle $\Omega_2$ in the 2D physical cell is mapped to its counterpart in the second quadrant.

Finally in Figure 7(d), by composition of operations in Figure 7(b) and Figure 7(c), $\square \mathrm{R}_4\mathrm{R}_5\mathrm{R}_6\mathrm{R}_7$ is mapped to $\square \mathrm{OR}_1\mathrm{R}_2\mathrm{R}_3$.

In summary, there should be four patches within $\square \mathrm{OR}_1\mathrm{R}_2\mathrm{R}_3$, namely, $(\mathrm{II} \cap \Omega_2) \cup \mathrm{I}, \mathrm{II} \cap \Omega_1, \mathrm{III} \cap \Omega_1, (\mathrm{III} \cap \Omega_2) \cup \mathrm{IV}$. Comparing Figure 7(a) with Figure 7(d), the linear mapping of each patch to $\square \mathrm{R}_4\mathrm{R}_5\mathrm{R}_6\mathrm{R}_7$ is $\mathcal{T}_0$, $\mathcal{T}_{-\mathbf{a}_1}$, $\mathcal{T}_{-\mathbf{a}_1-\mathbf{a}_2}$, $\mathcal{T}_{-\mathbf{a}_2}$, respectively.

Furthermore, comparing Figure 7(d), and Figure 2 we identify four patches Figure 7(d) with four patches within $\square \mathrm{OR}_1\mathrm{R}_2\mathrm{R}_3$ in Figure 2, namely

- $(\mathrm{II} \cap \Omega_2) \cup \mathrm{I} \mapsto \mathrm{I}, \quad \mathrm{II} \cap \Omega_1 \mapsto \mathrm{II},$
- $\mathrm{III} \cap \Omega_1 \mapsto \mathrm{III}, \quad (\mathrm{III} \cap \Omega_2) \cup \mathrm{IV} \mapsto \mathrm{IV}.$

## Appendix C. matrix $J_2$ and $J_3$ in triclinic lattice

### *Appendix C.1. Two cases of $J_2$*

As mentioned in (12), the BC (12) can be classified by the angle $\phi_3$ as Define

$$m_1 = \begin{cases} \text{floor}(\frac{a_2 \cos\phi_3}{\delta_x}), & \text{if } 0 < \phi_3 \le \frac{\pi}{2}, \\ \text{floor}(\frac{a_1 + a_2 \cos\phi_3}{\delta_x}), & \text{if } \frac{\pi}{2} < \phi_3 < \pi. \end{cases} \tag{C.1}$$

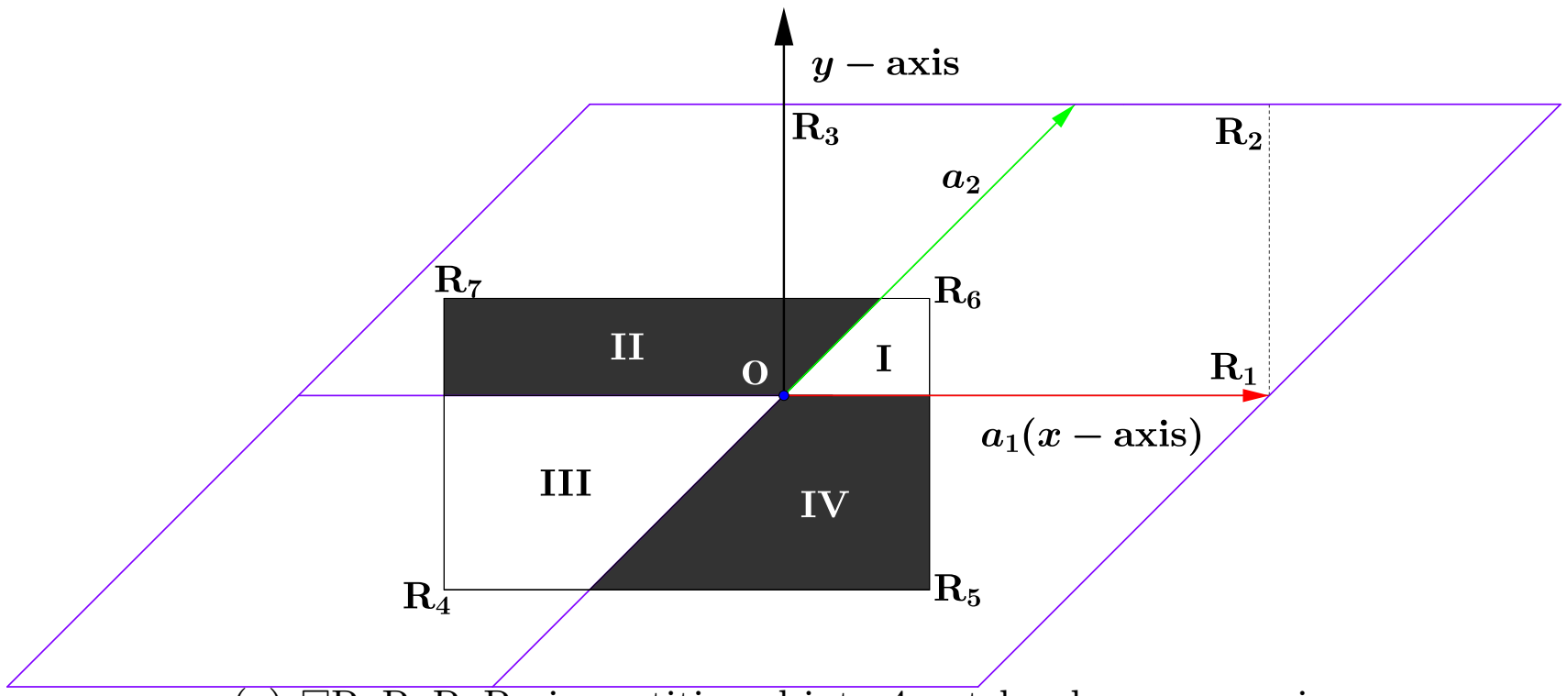

(a) $\square R_4R_5R_6R_7$ is partitioned into 4 patches by $\mathbf{a}_1$-,$\mathbf{a}_2$-axis.

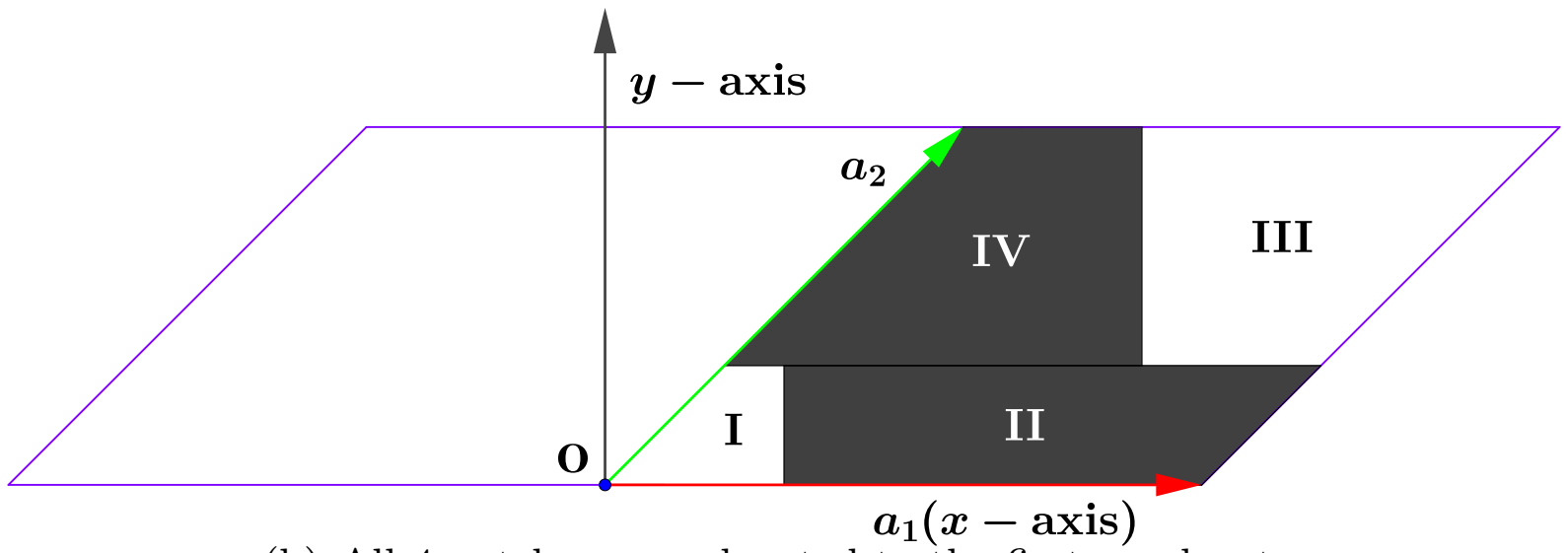

(b) All 4 patches are relocated to the first quadrant.

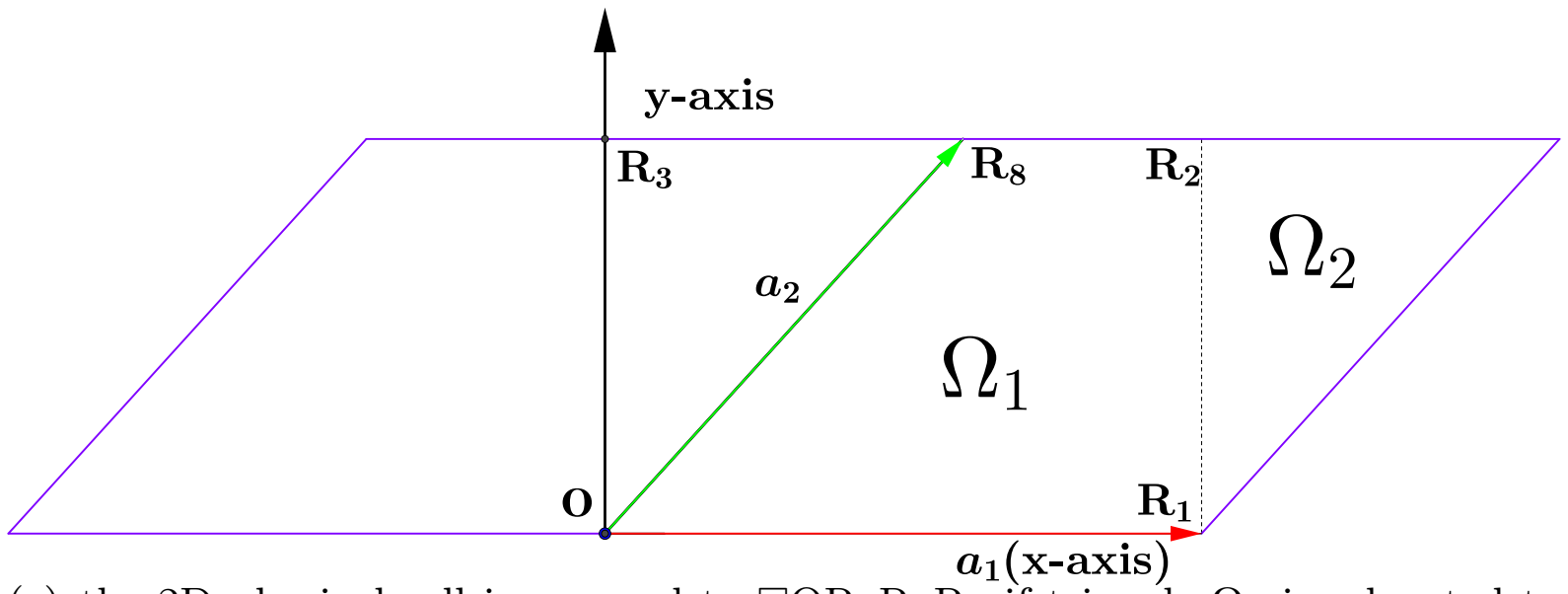

(c) the 2D physical cell is mapped to $\square OR_1R_2R_3$ if triangle $\Omega_2$ is relocated to triangle $OR_8R_3$.

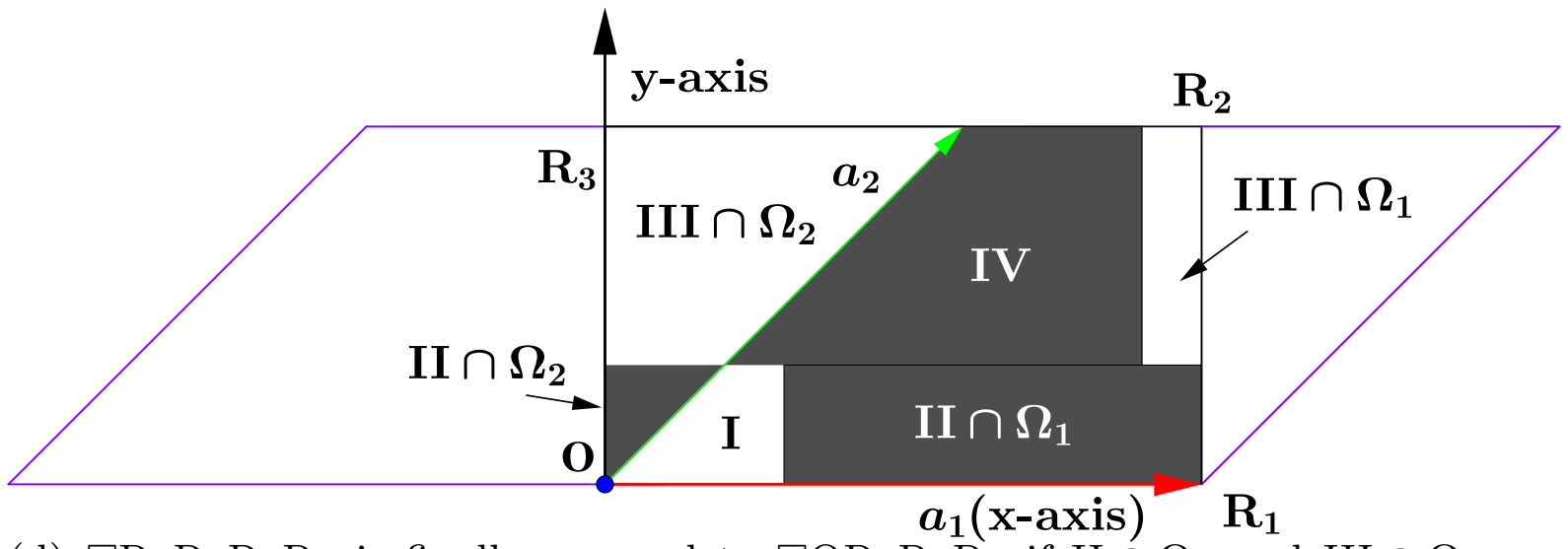

(d) $\square R_4R_5R_6R_7$ is finally mapped to $\square OR_1R_2R_3$ if $II \cap \Omega_2$ and $III \cap \Omega_2$ are relocated to the second quadrant.

Figure B.7: BC along $z$-direction.

Then $J_2$ in a triclinic system can be write as following cases

$$J_2 = \begin{cases} \xi(-\mathbf{k}\cdot\mathbf{a}_1)G_{n_1}(\mathbf{k}\cdot\mathbf{a}_1, n_1-m_1), & \text{if } 0<\phi_3\le\frac{\pi}{2}, \\ G_{n_1}(\mathbf{k}\cdot\mathbf{a}_1, n_1-m_1), & \text{if } \frac{\pi}{2}<\phi_3<\pi. \end{cases} \tag{C.2}$$

*Appendix C.2. Sixteen BCs in various lattice structure*

Recall that $\mathbf{a}_3^{\perp}$ is the projection of $\mathbf{a}_3$ onto $x$-$y$ plane. We classify the triclinic lattice into four categories according to the quadrant in which $\mathbf{a}_3^{\perp}$ is located, as shown in Figure C.8(1), C.9(2), C.10(3) and C.12(4). And according to the quadrant in which $\mathbf{a}_2$ is located and the first coordinate of $\mathbf{a}_1$, $\mathbf{a}_2$ $\mathbf{a}_3$, each category is further divided into four subcategories, as shown in Figure C.8(1-i), (1-ii), (1-iii), (1-iv), C.9(2-i), (2-ii), (2-iii), (2-iv), Figure C.10(3-i), (3-ii), (3-iii), (3-iv), Figure C.12(4-i), (4-ii), (4-iii) and (4-iv). Notice that the blue and green dotted vectors in these subfigures are equal to the translation vectors $\mathbf{a}_3$ and $\mathbf{a}_2$, respectively.

We first divide the top surface of $\mathbb{D}$ into red, green, and blue areas based on the categories described in Figure C.8(1), C.9(2), C.10(3) and C.12(4) It can be seen that the blue area has already fallen on the bottom of the working cell under $\mathcal{T}_{-\mathbf{a}_3}$, while the red and green areas want action of $\mathcal{T}_{\mathbf{a}_1}$ and/or $\mathcal{T}_{\mathbf{a}_2}$ in order to fall on the bottom surface of $\mathbb{D}$. The image of the top surface of $\mathbb{D}$ under $\mathcal{T}_{-\mathbf{a}_3}$ is partitioned into $\widetilde{\mathbf{I}}$, $\widetilde{\mathbf{II}}$, $\widetilde{\mathbf{III}}$, $\widetilde{\mathbf{IV}}$, while the bottom surface of $\mathbb{D}$ is partitioned into $\mathbf{I}$, $\mathbf{II}$, $\mathbf{III}$, $\mathbf{IV}$. We will discuss each subcategory, and reformulate the BC accordingly.

Let $\mathbf{x} = (x, y, 0) \in \mathbb{D}$ be the point in the bottom surface of $\mathbb{D}$.

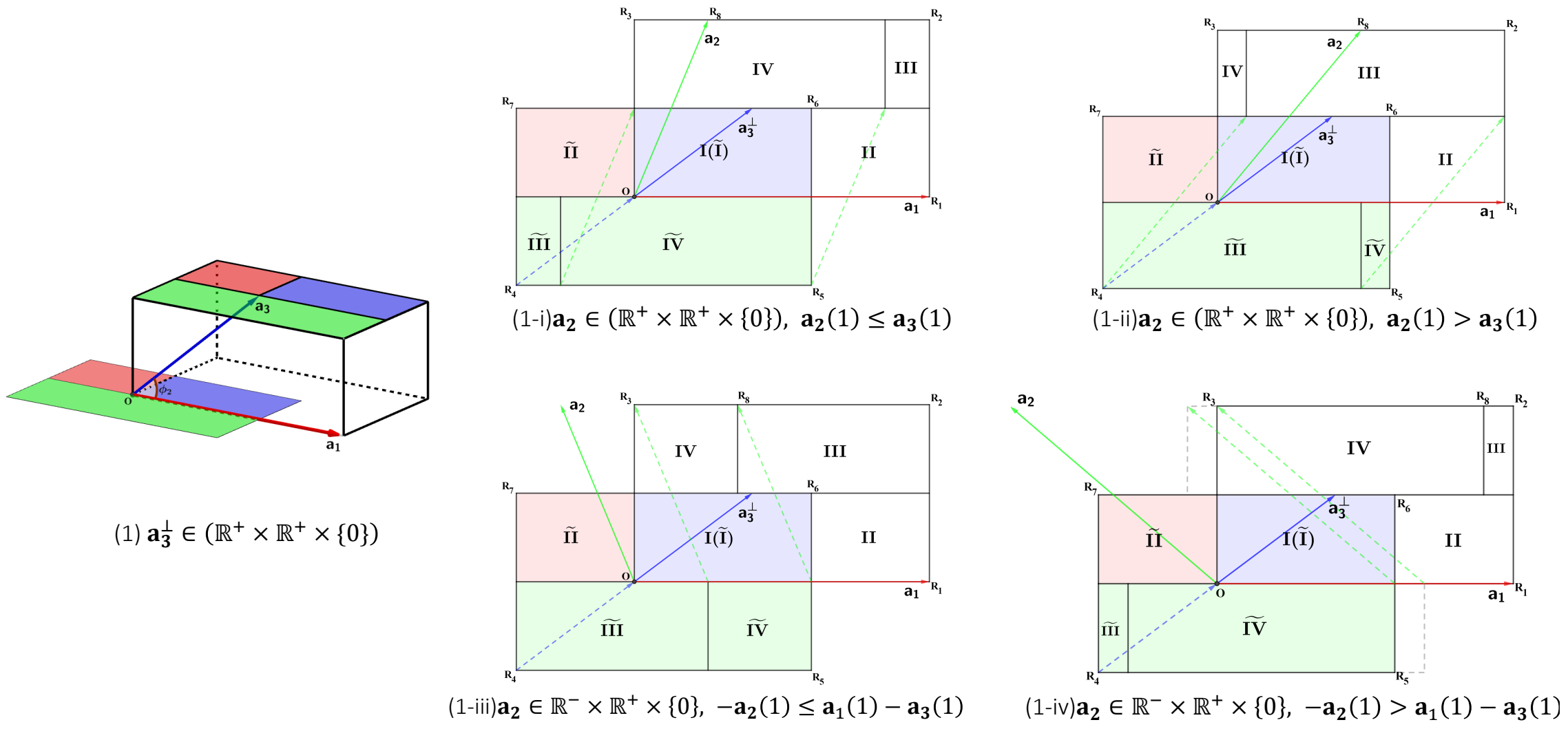

Figure C.8: Illustion of the first category in which $\mathbf{a}_3^{\perp}$ is located in the first quadrant.

- Case (1-i): $\mathbf{a}_3^\perp \in \mathbb{R}^+ \times \mathbb{R}^+ \times \{0\}$, $\mathbf{a}_2 \in \mathbb{R}^+ \times \mathbb{R}^+ \times \{0\}$, $\mathbf{a}_2(1) \le \mathbf{a}_3(1)$.

$$\boldsymbol{E}(\mathbf{x}) = \begin{cases} \boldsymbol{E}(\mathbf{x}), & \text{if } \mathbf{x} \in \mathbf{I} \\ \xi(\mathbf{k}\cdot\mathbf{a}_1)\boldsymbol{E}(\mathbf{x}-\mathbf{a}_1), & \text{if } \mathbf{x} \in \mathbf{II} \\ \xi(\mathbf{k}\cdot(\mathbf{a}_1+\mathbf{a}_2))\boldsymbol{E}(\mathbf{x}-\mathbf{a}_1-\mathbf{a}_2), & \text{if } \mathbf{x} \in \mathbf{III} \\ \xi(\mathbf{k}\cdot\mathbf{a}_2)\boldsymbol{E}(\mathbf{x}-\mathbf{a}_2), & \text{if } \mathbf{x} \in \mathbf{IV} \end{cases} \tag{C.3}$$

- Case (1-ii): $\mathbf{a}_3^\perp \in \mathbb{R}^+ \times \mathbb{R}^+ \times \{0\}$, $\mathbf{a}_2 \in \mathbb{R}^+ \times \mathbb{R}^+ \times \{0\}$, $\mathbf{a}_2(1) > \mathbf{a}_3(1)$

$$\boldsymbol{E}(\mathbf{x}) = \begin{cases} \boldsymbol{E}(\mathbf{x}), & \text{if } \mathbf{x} \in \mathbf{I} \\ \xi(\mathbf{k}\cdot\mathbf{a}_1)\boldsymbol{E}(\mathbf{x}-\mathbf{a}_1), & \text{if } \mathbf{x} \in \mathbf{II} \\ \xi(\mathbf{k}\cdot\mathbf{a}_2)\boldsymbol{E}(\mathbf{x}-\mathbf{a}_2), & \text{if } \mathbf{x} \in \mathbf{III} \\ \xi(\mathbf{k}\cdot(-\mathbf{a}_1+\mathbf{a}_2))\boldsymbol{E}(\mathbf{x}+\mathbf{a}_1-\mathbf{a}_2), & \text{if } \mathbf{x} \in \mathbf{IV} \end{cases} \tag{C.4}$$

- Case (1-iii): $\mathbf{a}_3^\perp \in \mathbb{R}^+ \times \mathbb{R}^+ \times \{0\}$, $\mathbf{a}_2 \in \mathbb{R}^- \times \mathbb{R}^+ \times \{0\}$, $-\mathbf{a}_2(1) \le \mathbf{a}_1(1) - \mathbf{a}_3(1)$.

$$\boldsymbol{E}(\mathbf{x}) = \begin{cases} \boldsymbol{E}(\mathbf{x}), & \text{if } \mathbf{x} \in \mathbf{I} \\ \xi(\mathbf{k}\cdot\mathbf{a}_1)\boldsymbol{E}(\mathbf{x}-\mathbf{a}_1), & \text{if } \mathbf{x} \in \mathbf{II} \\ \xi(\mathbf{k}\cdot(\mathbf{a}_1+\mathbf{a}_2))\boldsymbol{E}(\mathbf{x}-\mathbf{a}_1-\mathbf{a}_2), & \text{if } \mathbf{x} \in \mathbf{III} \\ \xi(\mathbf{k}\cdot\mathbf{a}_2)\boldsymbol{E}(\mathbf{x}-\mathbf{a}_2), & \text{if } \mathbf{x} \in \mathbf{IV} \end{cases} \tag{C.5}$$

- Case (1-iv): $\mathbf{a}_3^\perp \in \mathbb{R}^+ \times \mathbb{R}^+ \times \{0\}$, $\mathbf{a}_2 \in \mathbb{R}^- \times \mathbb{R}^+ \times \{0\}$, $-\mathbf{a}_2(1) > \mathbf{a}_1(1) - \mathbf{a}_3(1)$

$$\boldsymbol{E}(\mathbf{x}) = \begin{cases} \boldsymbol{E}(\mathbf{x}), & \text{if } \mathbf{x} \in \mathbf{I} \\ \xi(\mathbf{k}\cdot\mathbf{a}_1)\boldsymbol{E}(\mathbf{x}-\mathbf{a}_1), & \text{if } \mathbf{x} \in \mathbf{II} \\ \xi(\mathbf{k}\cdot(2\mathbf{a}_1+\mathbf{a}_2))\boldsymbol{E}(\mathbf{x}-2\mathbf{a}_1-\mathbf{a}_2), & \text{if } \mathbf{x} \in \mathbf{III} \\ \xi(\mathbf{k}\cdot(\mathbf{a}_1+\mathbf{a}_2))\boldsymbol{E}(\mathbf{x}-\mathbf{a}_1-\mathbf{a}_2), & \text{if } \mathbf{x} \in \mathbf{IV} \end{cases} \tag{C.6}$$

- Case (2-i): $\mathbf{a}_3^\perp \in \mathbb{R}^- \times \mathbb{R}^+ \times \{0\}$, $\mathbf{a}_2 \in \mathbb{R}^+ \times \mathbb{R}^+ \times \{0\}$, $\mathbf{a}_2(1) \le \mathbf{a}_1(1) + \mathbf{a}_3(1)$.

$$\boldsymbol{E}(\mathbf{x}) = \begin{cases} \xi(-\mathbf{k}\cdot\mathbf{a}_1)\boldsymbol{E}(\mathbf{x}+\mathbf{a}_1), & \text{if } \mathbf{x} \in \mathbf{I} \\ \boldsymbol{E}(\mathbf{x}), & \text{if } \mathbf{x} \in \mathbf{II} \\ \xi(\mathbf{k}\cdot\mathbf{a}_2)\boldsymbol{E}(\mathbf{x}-\mathbf{a}_2), & \text{if } \mathbf{x} \in \mathbf{III} \\ \xi(\mathbf{k}\cdot(-\mathbf{a}_1+\mathbf{a}_2))\boldsymbol{E}(\mathbf{x}+\mathbf{a}_1-\mathbf{a}_2), & \text{if } \mathbf{x} \in \mathbf{IV} \end{cases} \tag{C.7}$$

- Case (2-ii): $\mathbf{a}_3^\perp \in \mathbb{R}^- \times \mathbb{R}^+ \times \{0\}$, $\mathbf{a}_2 \in \mathbb{R}^+ \times \mathbb{R}^+ \times \{0\}$, $\mathbf{a}_2(1) > \mathbf{a}_1(1) + \mathbf{a}_3(1)$

$$\boldsymbol{E}(\mathbf{x}) = \begin{cases} \xi(-\mathbf{k}\cdot\mathbf{a}_1)\boldsymbol{E}(\mathbf{x}+\mathbf{a}_1), & \text{if } \mathbf{x} \in \mathbf{I} \\ \boldsymbol{E}(\mathbf{x}), & \text{if } \mathbf{x} \in \mathbf{II} \\ \xi(\mathbf{k}\cdot(-\mathbf{a}_1+\mathbf{a}_2))\boldsymbol{E}(\mathbf{x}+\mathbf{a}_1-\mathbf{a}_2), & \text{if } \mathbf{x} \in \mathbf{III} \\ \xi(\mathbf{k}\cdot(-2\mathbf{a}_1+\mathbf{a}_2))\boldsymbol{E}(\mathbf{x}+2\mathbf{a}_1-\mathbf{a}_2), & \text{if } \mathbf{x} \in \mathbf{IV} \end{cases} \tag{C.8}$$

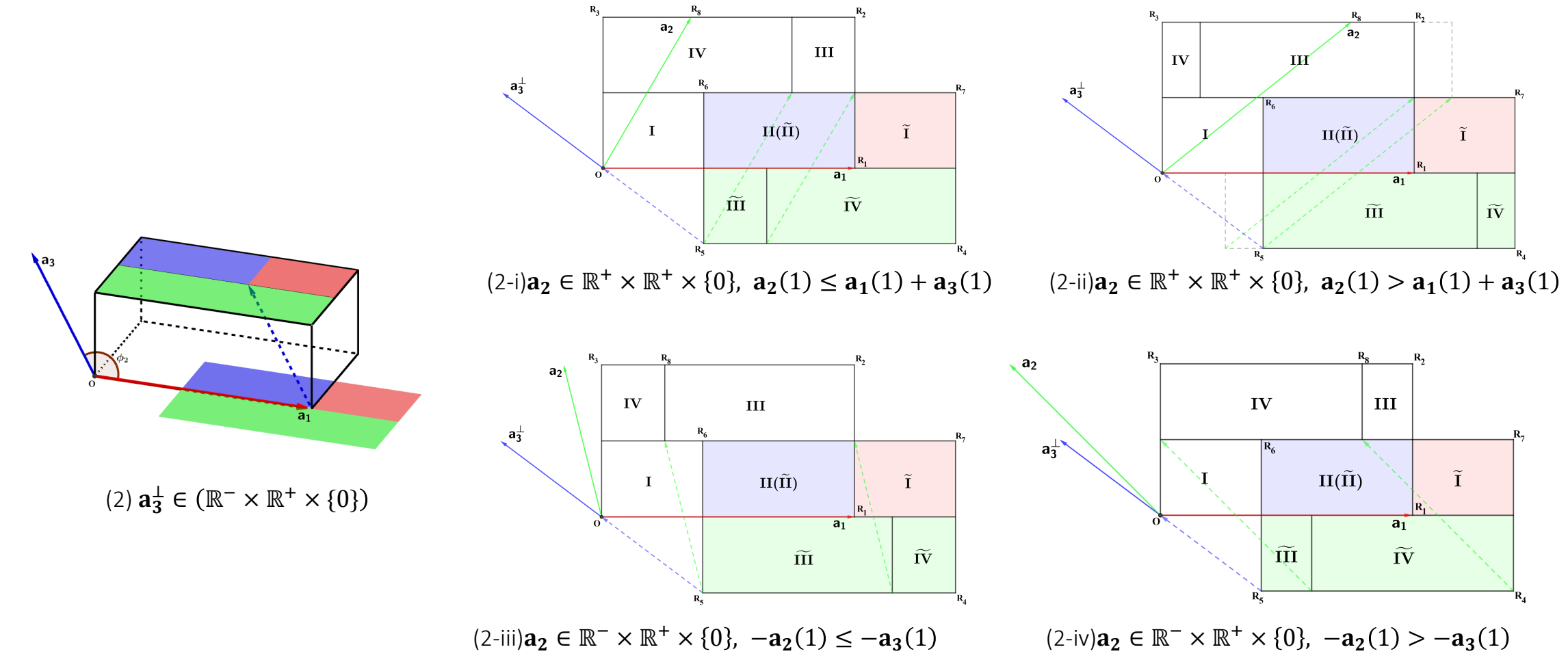

Figure C.9: Illustion of the second category in which $\mathbf{a}_3^\perp$ is located in the second quadrant.

- Case (2-iii): $\mathbf{a}_3^\perp \in \mathbb{R}^- \times \mathbb{R}^+ \times \{0\}$, $\mathbf{a}_2 \in \mathbb{R}^- \times \mathbb{R}^+ \times \{0\}$, $-\mathbf{a}_2(1) \le -\mathbf{a}_3(1)$.

$$\boldsymbol{E}(\mathbf{x}) = \begin{cases} \xi(-\mathbf{k}\cdot\mathbf{a}_1)\boldsymbol{E}(\mathbf{x}+\mathbf{a}_1), & \text{if } \mathbf{x} \in \mathbf{I} \\ \boldsymbol{E}(\mathbf{x}), & \text{if } \mathbf{x} \in \mathbf{II} \\ \xi(\mathbf{k}\cdot\mathbf{a}_2)\boldsymbol{E}(\mathbf{x}-\mathbf{a}_2), & \text{if } \mathbf{x} \in \mathbf{III} \\ \xi(\mathbf{k}\cdot(-\mathbf{a}_1+\mathbf{a}_2))\boldsymbol{E}(\mathbf{x}+\mathbf{a}_1-\mathbf{a}_2), & \text{if } \mathbf{x} \in \mathbf{IV} \end{cases} \tag{C.9}$$

- Case (2-iv): $\mathbf{a}_3^\perp \in \mathbb{R}^- \times \mathbb{R}^+ \times \{0\}$, $\mathbf{a}_2 \in \mathbb{R}^- \times \mathbb{R}^+ \times \{0\}$, $-\mathbf{a}_2(1) > -\mathbf{a}_3(1)$

$$\boldsymbol{E}(\mathbf{x}) = \begin{cases} \xi(-\mathbf{k}\cdot\mathbf{a}_1)\boldsymbol{E}(\mathbf{x}+\mathbf{a}_1), & \text{if } \mathbf{x} \in \mathbf{I} \\ \boldsymbol{E}(\mathbf{x}), & \text{if } \mathbf{x} \in \mathbf{II} \\ \xi(\mathbf{k}\cdot(\mathbf{a}_1+\mathbf{a}_2))\boldsymbol{E}(\mathbf{x}-\mathbf{a}_1-\mathbf{a}_2), & \text{if } \mathbf{x} \in \mathbf{III} \\ \xi(\mathbf{k}\cdot\mathbf{a}_2)\boldsymbol{E}(\mathbf{x}-\mathbf{a}_2), & \text{if } \mathbf{x} \in \mathbf{IV} \end{cases} \tag{C.10}$$

- Case (3-i): $\mathbf{a}_3^\perp \in \mathbb{R}^- \times \mathbb{R}^- \times \{0\}$, $\mathbf{a}_2 \in \mathbb{R}^+ \times \mathbb{R}^+ \times \{0\}$, $\mathbf{a}_2(1) \le -\mathbf{a}_3(1)$.

$$\boldsymbol{E}(\mathbf{x}) = \begin{cases} \xi(-\mathbf{k}\cdot(\mathbf{a}_1+\mathbf{a}_2))\boldsymbol{E}(\mathbf{x}+\mathbf{a}_1+\mathbf{a}_2), & \text{if } \mathbf{x} \in \mathbf{I} \\ \xi(-\mathbf{k}\cdot\mathbf{a}_2)\boldsymbol{E}(\mathbf{x}+\mathbf{a}_2), & \text{if } \mathbf{x} \in \mathbf{II} \\ \boldsymbol{E}(\mathbf{x}), & \text{if } \mathbf{x} \in \mathbf{III} \\ \xi(-\mathbf{k}\cdot\mathbf{a}_1)\boldsymbol{E}(\mathbf{x}+\mathbf{a}_1), & \text{if } \mathbf{x} \in \mathbf{IV} \end{cases} \tag{C.11}$$

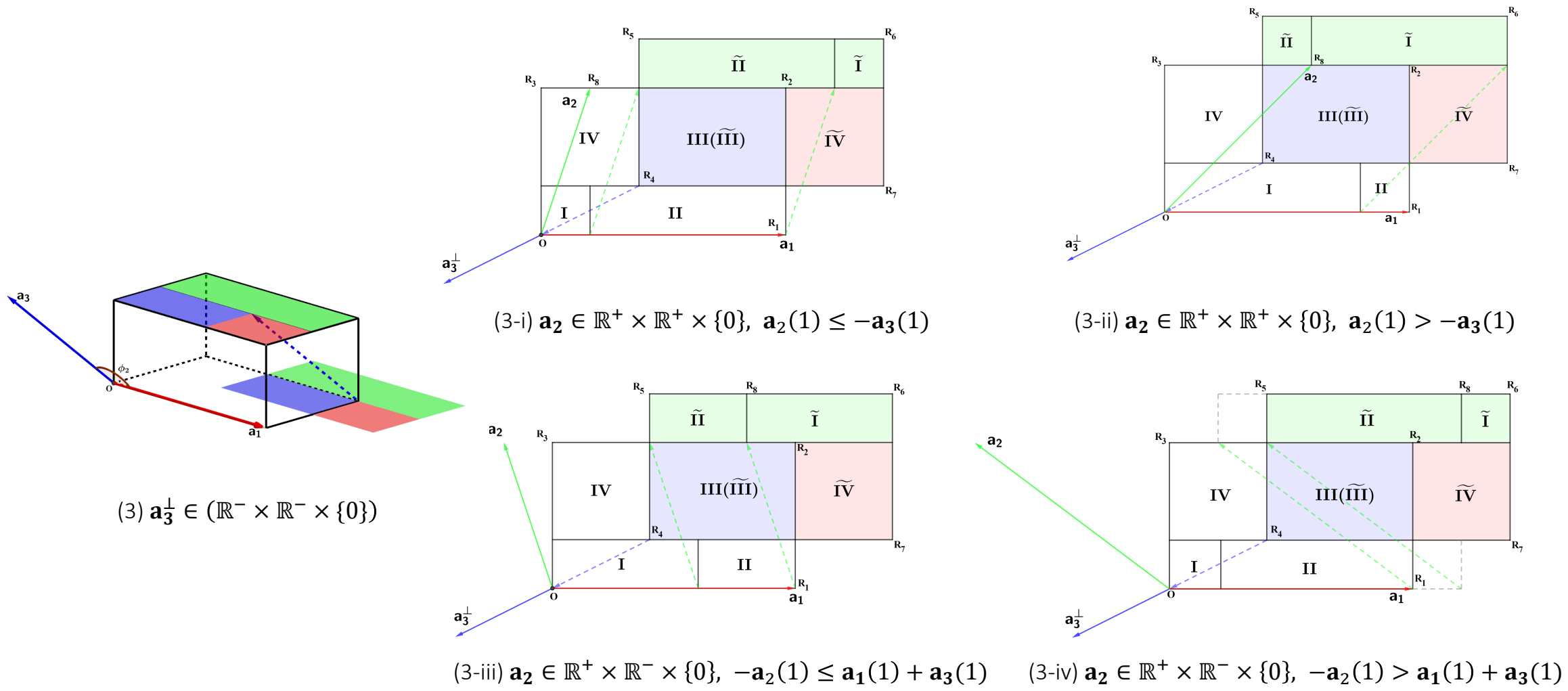

Figure C.10: Illustration of the third category in which $\mathbf{a}_3^\perp$ is located on the third quadrant.

- Case (3-ii): $\mathbf{a}_3^\perp \in \mathbb{R}^- \times \mathbb{R}^- \times \{0\}$, $\mathbf{a}_2 \in \mathbb{R}^+ \times \mathbb{R}^+ \times \{0\}$, $\mathbf{a}_2(1) > -\mathbf{a}_3(1)$

$$\boldsymbol{E}(\mathbf{x}) = \begin{cases} \xi(-\mathbf{k}\cdot\mathbf{a}_2)\boldsymbol{E}(\mathbf{x}+\mathbf{a}_2), & \text{if } \mathbf{x} \in \mathbf{I} \\ \xi(\mathbf{k}\cdot(\mathbf{a}_1-\mathbf{a}_2))\boldsymbol{E}(\mathbf{x}-\mathbf{a}_1+\mathbf{a}_2), & \text{if } \mathbf{x} \in \mathbf{II} \\ \boldsymbol{E}(\mathbf{x}), & \text{if } \mathbf{x} \in \mathbf{III} \\ \xi(-\mathbf{k}\cdot\mathbf{a}_1)\boldsymbol{E}(\mathbf{x}+\mathbf{a}_1), & \text{if } \mathbf{x} \in \mathbf{IV} \end{cases} \tag{C.12}$$

- Case (3-iii): $\mathbf{a}_3^\perp \in \mathbb{R}^- \times \mathbb{R}^- \times \{0\}$, $\mathbf{a}_2 \in \mathbb{R}^- \times \mathbb{R}^+ \times \{0\}$, $-\mathbf{a}_2(1) \le \mathbf{a}_1(1) + \mathbf{a}_3(1)$.

$$\boldsymbol{E}(\mathbf{x}) = \begin{cases} \xi(-\mathbf{k}\cdot(\mathbf{a}_1+\mathbf{a}_2))\boldsymbol{E}(\mathbf{x}+\mathbf{a}_1+\mathbf{a}_2), & \text{if } \mathbf{x} \in \mathbf{I} \\ \xi(-\mathbf{k}\cdot\mathbf{a}_2)\boldsymbol{E}(\mathbf{x}+\mathbf{a}_2), & \text{if } \mathbf{x} \in \mathbf{II} \\ \boldsymbol{E}(\mathbf{x}), & \text{if } \mathbf{x} \in \mathbf{III} \\ \xi(-\mathbf{k}\cdot\mathbf{a}_1)\boldsymbol{E}(\mathbf{x}+\mathbf{a}_1), & \text{if } \mathbf{x} \in \mathbf{IV} \end{cases} \tag{C.13}$$

- Case (3-iv): $\mathbf{a}_3^\perp \in \mathbb{R}^- \times \mathbb{R}^- \times \{0\}$, $\mathbf{a}_2 \in \mathbb{R}^- \times \mathbb{R}^+ \times \{0\}$, $-\mathbf{a}_2(1) > \mathbf{a}_1(1) + \mathbf{a}_3(1)$

$$\boldsymbol{E}(\mathbf{x}) = \begin{cases} \xi(-\mathbf{k}\cdot(2\mathbf{a}_1+\mathbf{a}_2))\boldsymbol{E}(\mathbf{x}+2\mathbf{a}_1+\mathbf{a}_2), & \text{if } \mathbf{x} \in \mathbf{I} \\ \xi(-\mathbf{k}\cdot(\mathbf{a}_1+\mathbf{a}_2))\boldsymbol{E}(\mathbf{x}+\mathbf{a}_1+\mathbf{a}_2), & \text{if } \mathbf{x} \in \mathbf{II} \\ \boldsymbol{E}(\mathbf{x}), & \text{if } \mathbf{x} \in \mathbf{III} \\ \xi(-\mathbf{k}\cdot\mathbf{a}_1)\boldsymbol{E}(\mathbf{x}+\mathbf{a}_1), & \text{if } \mathbf{x} \in \mathbf{IV} \end{cases} \tag{C.14}$$

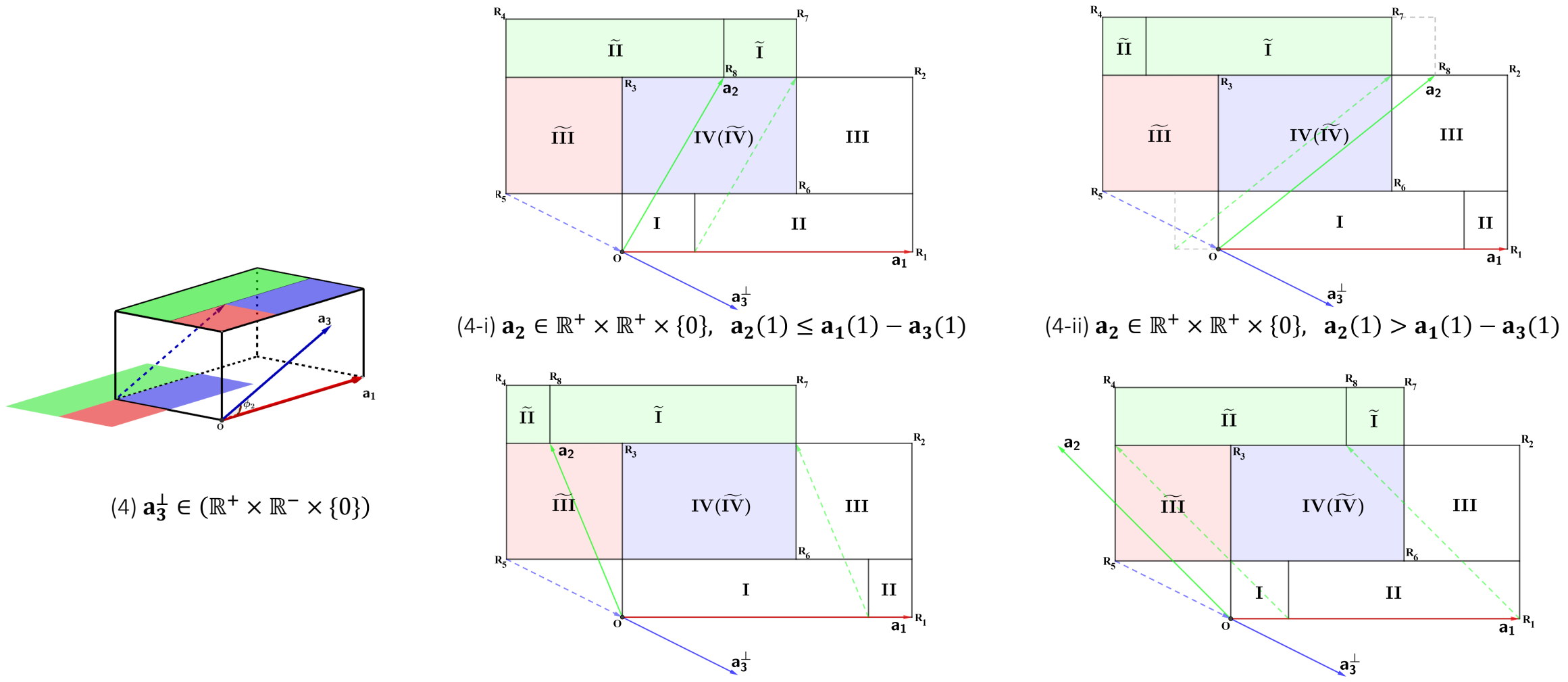

(4) $\mathbf{a}_3^{\perp} \in (\mathbb{R}^+ \times \mathbb{R}^- \times \{0\})$

(4-i) $\mathbf{a}_2 \in \mathbb{R}^+ \times \mathbb{R}^+ \times \{0\}$, $\mathbf{a}_2(1) \leq \mathbf{a}_1(1) - \mathbf{a}_3(1)$

(4-ii) $\mathbf{a}_2 \in \mathbb{R}^+ \times \mathbb{R}^+ \times \{0\}$, $\mathbf{a}_2(1) > \mathbf{a}_1(1) - \mathbf{a}_3(1)$

(4-iii) $\mathbf{a}_2 \in \mathbb{R}^+ \times \mathbb{R}^- \times \{0\}$, $-\mathbf{a}_2(1) \leq \mathbf{a}_3(1)$

(4-iv) $\mathbf{a}_2 \in \mathbb{R}^+ \times \mathbb{R}^- \times \{0\}$, $-\mathbf{a}_2(1) > \mathbf{a}_3(1)$

Figure C.11: Illustration of the fourth category in which $\mathbf{a}_3^{\perp}$ is located in the fourth quadrant.

- Case (4-i): $\mathbf{a}_3^{\perp} \in \mathbb{R}^+ \times \mathbb{R}^- \times \{0\}$, $\mathbf{a}_2 \in \mathbb{R}^+ \times \mathbb{R}^+ \times \{0\}$, $\mathbf{a}_2(1) \leq \mathbf{a}_1(1) - \mathbf{a}_3(1)$.

$$\boldsymbol{E}(\mathbf{x}) = \begin{cases} \xi(-\mathbf{k}\cdot\mathbf{a}_2)\boldsymbol{E}(\mathbf{x}+\mathbf{a}_2), & \text{if } \mathbf{x} \in \mathbf{I} \\ \xi(\mathbf{k}\cdot(\mathbf{a}_1-\mathbf{a}_2))\boldsymbol{E}(\mathbf{x}-\mathbf{a}_1+\mathbf{a}_2), & \text{if } \mathbf{x} \in \mathbf{II} \\ \xi(\mathbf{k}\cdot\mathbf{a}_1)\boldsymbol{E}(\mathbf{x}-\mathbf{a}_1), & \text{if } \mathbf{x} \in \mathbf{III} \\ \boldsymbol{E}(\mathbf{x}), & \text{if } \mathbf{x} \in \mathbf{IV} \end{cases} \tag{C.15}$$

- Case (4-ii): $\mathbf{a}_3^{\perp} \in \mathbb{R}^+ \times \mathbb{R}^- \times \{0\}$, $\mathbf{a}_2 \in \mathbb{R}^+ \times \mathbb{R}^+ \times \{0\}$, $\mathbf{a}_2(1) > \mathbf{a}_1(1) - \mathbf{a}_3(1)$

$$\boldsymbol{E}(\mathbf{x}) = \begin{cases} \xi(\mathbf{k}\cdot(\mathbf{a}_1-\mathbf{a}_2))\boldsymbol{E}(\mathbf{x}-\mathbf{a}_1+\mathbf{a}_2), & \text{if } \mathbf{x} \in \mathbf{I} \\ \xi(\mathbf{k}\cdot(2\mathbf{a}_1-\mathbf{a}_2))\boldsymbol{E}(\mathbf{x}-2\mathbf{a}_1+\mathbf{a}_2), & \text{if } \mathbf{x} \in \mathbf{II} \\ \xi(\mathbf{k}\cdot\mathbf{a}_1)\boldsymbol{E}(\mathbf{x}-\mathbf{a}_1), & \text{if } \mathbf{x} \in \mathbf{III} \\ \boldsymbol{E}(\mathbf{x}), & \text{if } \mathbf{x} \in \mathbf{IV} \end{cases} \tag{C.16}$$

- Case (4-iii): $\mathbf{a}_3^{\perp} \in \mathbb{R}^+ \times \mathbb{R}^- \times \{0\}$, $\mathbf{a}_2 \in \mathbb{R}^- \times \mathbb{R}^+ \times \{0\}$, $-\mathbf{a}_2(1) \leq \mathbf{a}_3(1)$.

$$\boldsymbol{E}(\mathbf{x}) = \begin{cases} \xi(-\mathbf{k}\cdot\mathbf{a}_2)\boldsymbol{E}(\mathbf{x}+\mathbf{a}_2), & \text{if } \mathbf{x} \in \mathbf{I} \\ \xi(\mathbf{k}\cdot(\mathbf{a}_1-\mathbf{a}_2))\boldsymbol{E}(\mathbf{x}-\mathbf{a}_1+\mathbf{a}_2), & \text{if } \mathbf{x} \in \mathbf{II} \\ \xi(\mathbf{k}\cdot\mathbf{a}_1)\boldsymbol{E}(\mathbf{x}-\mathbf{a}_1), & \text{if } \mathbf{x} \in \mathbf{III} \\ \boldsymbol{E}(\mathbf{x}), & \text{if } \mathbf{x} \in \mathbf{IV} \end{cases} \tag{C.17}$$

• Case (4-iv): $\mathbf{a}_3^\perp \in \mathbb{R}^+ \times \mathbb{R}^- \times \{0\}$, $\mathbf{a}_2 \in \mathbb{R}^- \times \mathbb{R}^+ \times \{0\}$, $-\mathbf{a}_2(1) > \mathbf{a}_3(1)$

$$\boldsymbol{E}(\mathbf{x}) = \begin{cases} \xi(-\mathbf{k}\cdot(\mathbf{a}_1+\mathbf{a}_2))\boldsymbol{E}(\mathbf{x}+\mathbf{a}_1+\mathbf{a}_2), & \text{if } \mathbf{x} \in \mathbf{I} \\ \xi(-\mathbf{k}\cdot\mathbf{a}_2)\boldsymbol{E}(\mathbf{x}+\mathbf{a}_2), & \text{if } \mathbf{x} \in \mathbf{II} \\ \xi(\mathbf{k}\cdot\mathbf{a}_1)\boldsymbol{E}(\mathbf{x}-\mathbf{a}_1), & \text{if } \mathbf{x} \in \mathbf{III} \\ \boldsymbol{E}(\mathbf{x}), & \text{if } \mathbf{x} \in \mathbf{IV} \end{cases} \tag{C.18}$$

In summary, the sixteen BCs (C.3-C.18) can be summarized into the following equation:

$$\boldsymbol{E}(\mathbf{x}) = \begin{cases} \xi(-\mathbf{k}\cdot\mathbf{t}_1)\boldsymbol{E}(\mathbf{x}+\mathbf{t}_1), & \text{if } \mathbf{x} \in \mathbf{I} \\ \xi(-\mathbf{k}\cdot\mathbf{t}_2)\boldsymbol{E}(\mathbf{x}+\mathbf{t}_2), & \text{if } \mathbf{x} \in \mathbf{II} \\ \xi(-\mathbf{k}\cdot\mathbf{t}_3)\boldsymbol{E}(\mathbf{x}+\mathbf{t}_3), & \text{if } \mathbf{x} \in \mathbf{III} \\ \xi(-\mathbf{k}\cdot\mathbf{t}_4)\boldsymbol{E}(\mathbf{x}+\mathbf{t}_4), & \text{if } \mathbf{x} \in \mathbf{IV} \end{cases} \tag{C.19}$$

where $\{\mathbf{t}_i\}_{i=1}^4$ can be substituted by the translation vectors in (C.3-C.18).

*Appendix C.3. General formula of $J_3$*

Define

$$m_2 = \text{floor}(\frac{\overline{R_5R_1}}{\delta_x}),\ m_3 = \text{floor}(\frac{\overline{R_4R_3}}{\delta_y}),\ m_4 = \text{floor}(\frac{\overline{R_6R_2}}{\delta_x}), \tag{C.20}$$

then the matrix $J_3$ corresponding to BCs (C.3-C.18) can be expressed by the following formulation

$$J_3 = \begin{bmatrix} 0 & I_{m_3} \otimes \begin{bmatrix} 0 & \xi(\mathbf{k}\cdot\mathbf{t}_3)I_{m_4} \\ \xi(\mathbf{k}\cdot\mathbf{t}_4)I_{n_1-m_4} & 0 \end{bmatrix} \\ I_{n_2-m_3} \otimes \begin{bmatrix} 0 & \xi(\mathbf{k}\cdot\mathbf{t}_2)I_{m_2} \\ \xi(\mathbf{k}\cdot\mathbf{t}_1)I_{n_1-m_2} & 0 \end{bmatrix} & 0 \end{bmatrix} \tag{C.21}$$

## Appendix D. matrix $J_3$ for all Bravais lattices except triclinic lattice

Matrix $J_3$ can be represented as the general form

$$J_3 = \begin{bmatrix} 0 & \eta_1 I_{m_3} \otimes \begin{bmatrix} 0 & \eta_2 I_{k_1} \\ \eta_3 I_{n_1-k_1} & 0 \end{bmatrix} \\ \eta_4 I_{n_2-m_3} \otimes \begin{bmatrix} 0 & \eta_5 I_{k_2} \\ \eta_6 I_{n_1-k_2} & 0 \end{bmatrix} & 0 \end{bmatrix} \tag{D.1}$$

with unimodular $\eta_i$ for $i = 1, \ldots, 6$ and $k_1, k_2 \in \mathbb{N}$. We provide specific expressions of $\eta_i, i = 1, \ldots, 6$ and $k_1, k_2$ for all Bravais lattices except triclinic lattice below. Denote $\zeta_1 = \exp(-\imath 2\pi\mathbf{k}\cdot\mathbf{a}_1)$, $\zeta_2 = \exp(-\imath 2\pi\mathbf{k}\cdot\mathbf{a}_2)$, $\zeta_3 = \exp(\imath 2\pi\mathbf{k}\cdot\mathbf{a}_3)$. Note that any $\eta_i$ that is not specified below is just 1.

• Cubic system

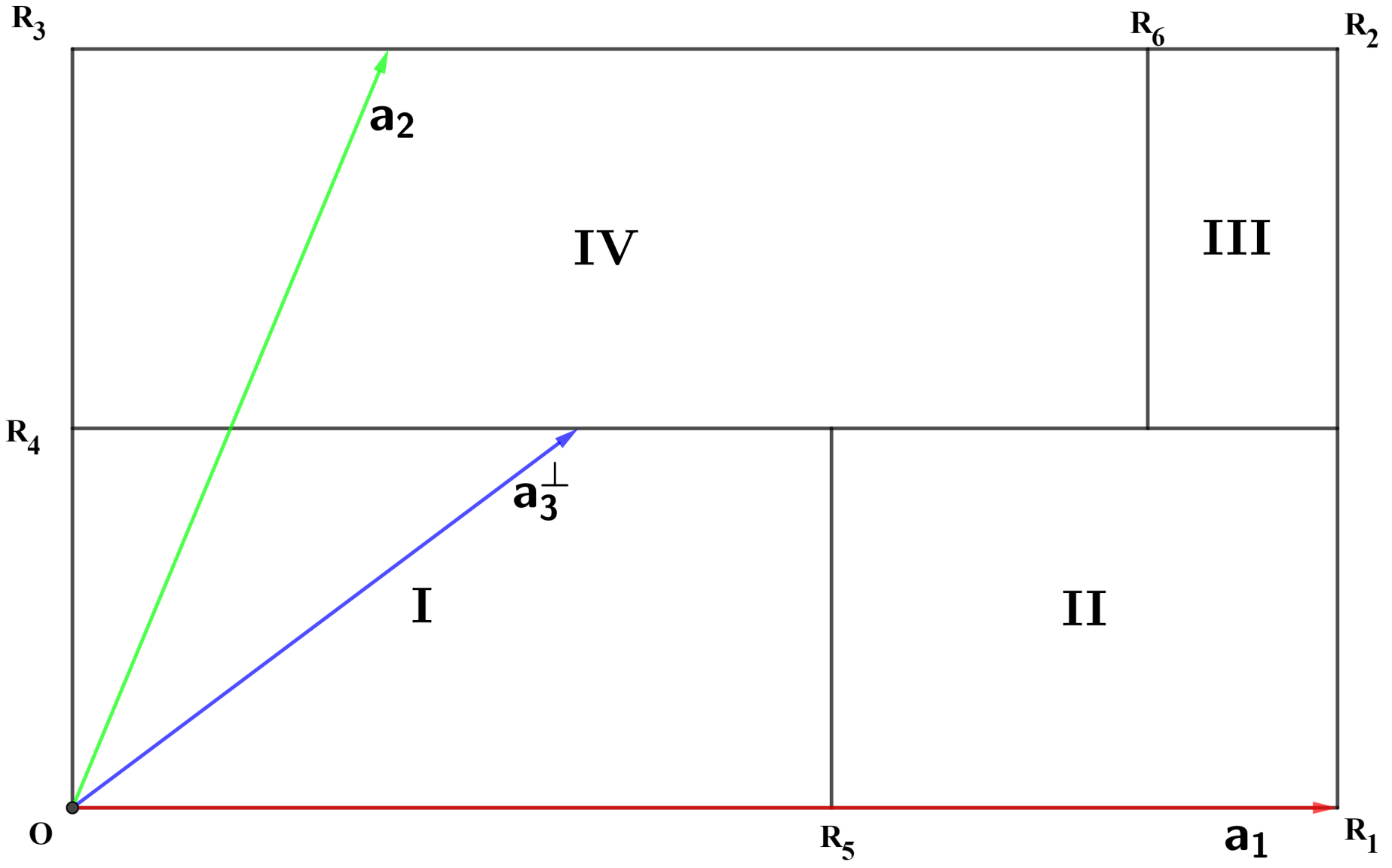

Figure C.12: The bottom surface of working cell in subcategory (1-i).

**(1)** Primitive: $J_3 = I_{n_1 n_2}$.

**(2)** F.C. (Here $\phi_3 < \pi/2$, $\phi_2 < \pi/2$, $\ell_2 > 0$, $m_1 = m_2 = n_1/2$, $m_3 = n_2/3$):

$$\eta_1 = \zeta_2,\ \eta_5 = \zeta_1, \quad k_1 = 0,\ k_2 = m_2.$$

**(3)** B.C. (Here $\phi_3 > \pi/2$, $\phi_2 > \pi/2$, $\ell_2 < 0$, $m_1 = m_2 = 2n_1/3$, $m_3 = n_2/2$):

$$\eta_3 = \eta_6 = \zeta_1^{-1},\ \eta_4 = \zeta_2^{-1}, \quad k_1 = m_2,\ k_2 = n_1 - m_2.$$

• Hexagonal system

**(a)** Assuming $a_1 \geq a_3$ (Here $\phi_3 > \pi/2$, $\phi_2 = \pi/2$, $\phi_1 = \pi/2$, $m_1 = n_1/2$, $m_2 = m_3 = 0$): $J_3 = I_{n_1 n_2}$.

**(b)** Assuming $a_1 < a_3$ (Here $\phi_3 = \pi/2,\ \phi_2 = \pi/2,\ \phi_1 > \pi/2,\ m_1 = m_2 = 0,$ $m_3 = n_2/2$):

$$\eta_4 = \zeta_2^{-1}, \quad k_1 = k_2 = 0.$$

• Rhombohedral system

**(a)** Assuming $\sqrt{2}a_3 < \sqrt{3}a_1$ (Here $\phi_3 < \pi/2,\ \phi_2 < \pi/2,\ \ell_2 > 0,\ m_1 = m_2 \geq n_1/2$):

$$\eta_1 = \zeta_2,\ \eta_5 = \zeta_1, \quad k_1 = 0,\ k_2 = m_2.$$

**(b)** Assuming $\sqrt{2}a_3 > \sqrt{3}a_1$ (Here $\phi_3 > \pi/2,\ \phi_2 > \pi/2,\ \ell_2 < 0,\ m_1 = m_2 \geq n_1/2$):

$$\eta_3 = \eta_6 = \zeta_1^{-1},\ \eta_4 = \zeta_2^{-1}, \quad k_1 = m_2,\ k_2 = m_1 + m_2 - n_1.$$

• Tetragonal system

**(1)** Primitive: $J_3 = I_{n_1 n_2}$.

**(2)** Body-Centered:

**(a)** Assuming $a_3 \leq \sqrt{2}a_1$ (Here $\phi_3 > \pi/2$, $\phi_2 > \pi/2$, $\ell_2 \leq 0$, $m_1 = m_2$):

$$\eta_3 = \eta_6 = \zeta_1^{-1},\ \eta_4 = \zeta_2^{-1}, \quad k_1 = m_2,\ k_2 = m_1 + m_2 - n_1.$$

**(b)** Assuming $a_3 > \sqrt{2}a_1$ (Here $\phi_3 > \pi/2$, $\phi_2 > \pi/2$, $\ell_2 > 0$, $m_1 = m_2$):

$$\eta_3 = \eta_6 = \zeta_1^{-1},\ \eta_4 = \zeta_1^{-1}\zeta_2^{-1}, \quad k_1 = m_2,\ k_2 = m_1 + m_2.$$

• Orthorhombic system

**(1)** Primitive: $J_3 = I_{n_1 n_2}$.

**(2)** A-Base-centered (Here $\phi_3 = \pi/2$, $\phi_2 = \pi/2$, $\phi_1 < \pi/2$, $m_1 = m_2 = 0$):

$$\eta_1 = \zeta_2, \quad k_1 = k_2 = 0.$$

**(3)** C-Base-centered (Here $\phi_3 > \pi/2$, $\phi_2 = \pi/2$, $\phi_1 = \pi/2$, $m_1 = m_2 = 0$):

$$J_3 = I_{n_1 n_2}.$$

**(4)** Face-Centered (Here $\phi_3 < \pi/2$, $\phi_2 < \pi/2$, $\ell_2 > 0$, $m_1 > m_2$):

$$\eta_1 = \zeta_2,\ \eta_3 = \zeta_1^{-1},\ \eta_5 = \zeta_1, \quad k_1 = n_1 - m_1 + m_2,\ k_2 = m_2.$$

**(5)** Body-Centered

**(a)** Assuming $a_1 \geq \sqrt{a_2^2 + a_3^2}$ (Here $\phi_3 > \pi/2,\ \phi_2 > \pi/2,\ \ell_2 \geq 0,\ m_1 < m_2$):

$$\eta_1 = \zeta_2,\ \eta_2 = \zeta_1,\ \eta_6 = \zeta_1^{-1}, \quad k_1 = m_2 - m_1,\ k_2 = m_2.$$

**(b)** Assuming $a_1 < \sqrt{a_2^2 + a_3^2}$ (Here $\phi_3 > \pi/2$, $\phi_2 > \pi/2$, $\ell_2 < 0$, $m_1 + m_2 > n_1$):

$$\eta_3 = \eta_6 = \zeta_1^{-1},\ \eta_4 = \zeta_2^{-1}, \quad k_1 = m_2,\ k_2 = m_1 + m_2 - n_1.$$

- Monoclinic system

**(1)** Primitive

**(a)** Assuming $a_1 \geq a_3$ (Here $\phi_2 = \phi_1 = \pi/2$, $m_2 = m_3 = 0$): $J_3 = I_{n_1 n_2}$.

**(b)** Assuming $a_1 < a_3$ and $\phi_3 < \pi/2$ (Here $\phi_3 = \phi_2 = \pi/2$, $\phi_1 < \pi/2$, $m_2 = m_3 = 0$):

$$\eta_1 = \zeta_2, \quad k_1 = k_2 = 0.$$

**(c)** Assuming $a_1 < a_3$ and $\phi_3 > \pi/2$ (Here $\phi_3 = \pi/2$, $\phi_2 = \pi/2$, $\phi_1 > \pi/2$, $m_1 = m_2 = 0$):

$$\eta_4 = \zeta_2^{-1}, \quad k_1 = k_2 = 0.$$

**(2)** Base-Centered

**(a)** Assuming $a_1 \geq \sqrt{a_2^2 + a_3^2}/2$, $\phi_3 > \pi/2$, $a_2 \sin\phi_3 \geq a_3$ (Here $\phi_3 < \pi/2$, $\phi_2 < \pi/2$, $\ell_2 \geq 0$, $m_1 = m_2$):

$$\eta_1 = \zeta_2,\ \eta_5 = \zeta_1, \quad k_1 = 0,\ k_2 = m_2.$$

**(b)** Assuming $a_1 \geq \sqrt{a_2^2 + a_3^2}/2$, $\phi_3 > \pi/2$, $a_2 \sin\phi_3 < a_3$ (Here $\phi_3 < \pi/2$, $\phi_2 < \pi/2$, $\ell_2 < 0$, $m_1 + m_2 < n_1$):

$$\eta_2 = \eta_5 = \zeta_1,\ \eta_4 = \zeta_2^{-1}, \quad k_1 = m_2,\ k_2 = m_1 + m_2.$$

**(c)** Assuming $a_1 \geq \sqrt{a_2^2 + a_3^2}/2$, $\phi_3 < \pi/2$, $a_2 \sin\phi_3 \geq a_3$ (Here $\phi_3 > \pi/2$, $\phi_2 > \pi/2$, $\cos\phi_1 \geq \cos\phi_3 \cos\phi_2$, $m_1 = m_2$):

$$\eta_1 = \zeta_2,\ \eta_6 = \zeta_1^{-1}, \quad k_1 = 0,\ k_2 = m_2.$$

**(d)** Assuming $a_1 \geq \sqrt{a_2^2 + a_3^2}/2$, $\phi_3 < \pi/2$, $a_2 \sin\phi_3 < a_3$ (Here $\phi_3 > \pi/2$, $\phi_2 > \pi/2$, $\ell_2 < 0$, $m_1 + m_2 > n_1$):

$$\eta_3 = \eta_6 = \zeta_1^{-1},\ \eta_4 = \zeta_2^{-1}, \quad k_1 = m_2,\ k_2 = m_1 + m_2 - n_1.$$

**(e)** Assuming $a_1 \geq \sqrt{a_2^2 + a_3^2}/2$, $\phi_3 > \pi/2$, $a_2 \geq a_3$ (Here $\phi_3 \leq \pi/2$, $\phi_2 < \pi/2$, $\ell_2 > 0$):
If $m_1 \leq m_2$, then

$$\eta_1 = \zeta_2,\ \eta_2 = \eta_5 = \zeta_1, \quad k_1 = m_2 - m_1,\ k_2 = m_2;$$

otherwise,

$$\eta_1 = \zeta_2,\ \eta_3 = \zeta_1^{-1},\ \eta_5 = \zeta_1, \quad k_1 = n_1 - m_1 + m_2,\ k_2 = m_2.$$

**(f)** Assuming $a_1 < \sqrt{a_2^2 + a_3^2}/2,\ \phi_3 > \pi/2,\ a_2 < a_3$ (Here $\phi_3 > \pi/2$, $\phi_2 < \pi/2,\ \ell_2 > 0$):
If $m_1 \le m_2$, then

$$\eta_1 = \zeta_1\zeta_2,\ \eta_2 = \eta_5 = \zeta_1, \quad k_1 = m_2 - m_1,\ k_2 = m_2;$$

otherwise,

$$\eta_1 = \zeta_2,\ \eta_2 = \eta_5 = \zeta_1, \quad k_1 = n_1 - m_1 + m_2,\ k_2 = m_2.$$

**(g)** Assuming $a_1 < \sqrt{a_2^2 + a_3^2}/2,\ \phi_3 < \pi/2,\ a_2 \ge a_3$ (Here $\phi_3 \le \pi/2$, $\phi_2 > \pi/2,\ \ell_2 < 0$):
If $m_1 + m_2 \le n_1$, then

$$\eta_3 = \eta_6 = \zeta_1^{-1},\ \eta_4 = \zeta_2^{-1}, \quad k_1 = m_2,\ k_2 = m_1 + m_2;$$

otherwise,

$$\eta_3 = \zeta_1^{-1},\ \eta_4 = \zeta_2^{-1},\ \eta_5 = \zeta_1, \quad k_1 = m_2,\ k_2 = m_1 + m_2 - n_1.$$

**(h)** Assuming $a_1 \ge \sqrt{a_2^2 + a_3^2}/2,\ \phi_3 < \pi/2,\ a_2 < a_3$ (Here $\phi_3 > \pi/2$, $\phi_2 > \pi/2,\ \ell_2 < 0$):
If $m_1 + m_2 \le n_1$, then

$$\eta_3 = \eta_6 = \zeta_1^{-1},\ \eta_4 = \zeta_1^{-1}\zeta_2^{-1}, \quad k_1 = m_2,\ k_2 = m_1 + m_2;$$

otherwise,

$$\eta_3 = \eta_6 = \zeta_1^{-1},\ \eta_4 = \zeta_2^{-1}, \quad k_1 = m_2,\ k_2 = m_1 + m_2 - n_1.$$

## Appendix E. Fast algorithms for $T\mathbf{q}$ and $T^*\mathbf{p}$

$T\mathbf{q}$ and $T^*\mathbf{p}$ are cornerstone of our fast eigensolver in Figure 4, which can be computed in $\mathcal{O}(n \log n)$ flops using Algorithms 1 and 2.

Here are the definitions of some matrices in Algorithms 1 and 2

$$\begin{aligned}
[U_m]_{ij} &= \xi((i-1)j/m),\ i,j = 1,2,\cdots,m,\\
[D_{\mathbf{a}_1}]_{ij} &= \xi((i-1)\theta_{\mathbf{a}_1})\delta_{ij},\ i,j = 1,2,\cdots,n_1,\\
[D_{\hat{\mathbf{a}}_2}]_{ji} &= \xi((j-1)\theta_{\hat{\mathbf{a}}_2,i}),\ i = 1,2,\cdots,n_1,\ j = 1,2,\cdots,n_2,\\
[D_{\hat{\mathbf{a}}_3,i}]_{kj} &= \xi((k-1)\theta_{\hat{\mathbf{a}}_3,ij}),\ i = 1,2,\cdots,n_1,\ j = 1,2,\cdots,n_2,\ k = 1,2,\cdots,n_3,
\end{aligned}$$

and $\circ$ refers to the Hadamard product, *i.e.*, pointwise product.

**Algorithm 1** FFT-based matrix-vector multiplication for $T\mathbf{q}$

**Input:** Any vector $\mathbf{q} = \begin{bmatrix}\mathbf{q}_1^\top & \cdots & \mathbf{q}_{n_1}^\top\end{bmatrix}^\top \in \mathbb{C}^n$ with $\mathbf{q}_i = \begin{bmatrix}\mathbf{q}_{i,1}^\top & \cdots & \mathbf{q}_{i,n_2}^\top\end{bmatrix}^\top$ and $\mathbf{q}_{ij} \in \mathbb{C}^{n_3}$ for $i = 1, \ldots, n_1, j = 1, \ldots, n_2$.

**Output:** The vector $\mathbf{g} \equiv T\mathbf{q}$.

Set $\widetilde{Q}_{\mathbf{z},i} = \begin{bmatrix}\mathbf{q}_{i,1} & \cdots & \mathbf{q}_{i,n_2}\end{bmatrix}$ and $\widetilde{Q}_{\mathbf{z}} = \begin{bmatrix}\widetilde{Q}_{\mathbf{z},1} & \cdots & \widetilde{Q}_{\mathbf{z},n_1}\end{bmatrix}$;

Compute $\widetilde{Q}_{u\mathbf{z}} = U_{n_3}\widetilde{Q}_{\mathbf{z}}$ by backward FFT;

Compute $\widetilde{Q}_{\mathbf{z},i} = D_{\hat{\mathbf{a}}_3,i} \circ \widetilde{Q}_{u\mathbf{z}}(:, (i-1)n_2 + 1 : in_2)$ for $i = 1, \ldots, n_1$;

Set $\widetilde{Q}_{\mathbf{y}} = \begin{bmatrix}\widetilde{Q}_{\mathbf{z},1}^\top & \cdots & \widetilde{Q}_{\mathbf{z},n_1}^\top\end{bmatrix} \in \mathbb{C}^{n_2 \times n_1 n_3}$;

Compute $\widetilde{Q}_{u\mathbf{y}} \equiv \begin{bmatrix}\widetilde{Q}_{u\mathbf{y}}^{(1)} & \cdots & \widetilde{Q}_{u\mathbf{y}}^{(n_1)}\end{bmatrix} = U_{n_2}\widetilde{Q}_{\mathbf{y}}$ by backward FFT;

Compute $\widetilde{Q}_{\mathbf{y},k} = D_{\hat{\mathbf{a}}_2} \circ \begin{bmatrix}\widetilde{Q}_{u\mathbf{y}}^{(1)}(:, k) & \cdots & \widetilde{Q}_{u\mathbf{y}}^{(n_1)}(:, k)\end{bmatrix}$ for $k = 1, \ldots, n_3$;

Set $\widetilde{Q}_{\mathbf{x}} = \begin{bmatrix}\widetilde{Q}_{\mathbf{y},1}^\top & \cdots & \widetilde{Q}_{\mathbf{y},n_3}^\top\end{bmatrix}$;

Compute $\widetilde{Q}_{u\mathbf{x}} = U_{n_1}\widetilde{Q}_{\mathbf{x}}$ by backward FFT;

Compute $\mathbf{g} = D_{\mathbf{a}_1}\widetilde{Q}_{u\mathbf{x}}/\sqrt{n_1 n_2 n_3}$; $\quad \mathbf{g} = \mathbf{g}(:)$.

**Algorithm 2** FFT-based matrix-vector multiplication for $T^*\mathbf{p}$

**Input:** Any vector $\mathbf{p} = \begin{bmatrix}\mathbf{p}_1^\top & \cdots & \mathbf{p}_{n_3}^\top\end{bmatrix}^\top \in \mathbb{C}^n$ with $\mathbf{p}_k = \begin{bmatrix}\mathbf{p}_{1,k}^\top & \cdots & \mathbf{p}_{n_2,k}^\top\end{bmatrix}^\top$ and $\mathbf{p}_{j,k} \in \mathbb{C}^{n_1}$ for $j = 1, \ldots, n_2$, $k = 1, \ldots, n_3$.

**Output:** The vector $\mathbf{f} \equiv T^*\mathbf{p}$.

Set $\widetilde{P}_{\mathbf{x},k} = \begin{bmatrix}\mathbf{p}_{1,k} & \cdots & \mathbf{p}_{n_2,k}\end{bmatrix}$ and $\widetilde{P}_{\mathbf{x}} = \begin{bmatrix}\widetilde{P}_{\mathbf{x},1} & \cdots & \widetilde{P}_{\mathbf{x},n_3}\end{bmatrix}$;

Compute $\widetilde{P}_{e\mathbf{x}} = D_{\mathbf{a}_1}^*\widetilde{P}_{\mathbf{x}}$;

Compute $\widetilde{P}_{u\mathbf{x}} = U_{n_1}^*\widetilde{P}_{e\mathbf{x}} \in \mathbb{C}^{n_1 \times n_2 n_3}$ by forward FFT;

Set $\widetilde{P}_{u\mathbf{x}}^{(i)} = \begin{bmatrix}\widetilde{P}_{u\mathbf{x}}(i, 1 : n_2)^\top & \cdots & \widetilde{P}_{u\mathbf{x}}(i, (n_3 - 1)n_2 + 1 : n_3 n_2)^\top\end{bmatrix}$;

Compute $\widetilde{P}_{e\mathbf{y}} = \begin{bmatrix}D_{\hat{\mathbf{a}}_2,1}^*\widetilde{P}_{u\mathbf{x}}^{(1)} & \cdots & D_{\hat{\mathbf{a}}_2,n_1}^*\widetilde{P}_{u\mathbf{x}}^{(n_1)}\end{bmatrix}$;

Compute $\widetilde{P}_{u\mathbf{y}} \equiv \begin{bmatrix}\widetilde{P}_{u\mathbf{y}}^{(1)} & \cdots & \widetilde{P}_{u\mathbf{y}}^{(n_1)}\end{bmatrix} = U_{n_2}^*\widetilde{P}_{e\mathbf{y}} \in \mathbb{C}^{n_2 \times n_1 n_3}$ by forward FFT;

Compute $\widetilde{P}_{e\mathbf{z}} = \begin{bmatrix}\bar{D}_{\hat{\mathbf{a}}_3,1} \circ (\widetilde{P}_{u\mathbf{y}}^{(1)})^\top & \cdots & \bar{D}_{\hat{\mathbf{a}}_3,n_1} \circ (\widetilde{P}_{u\mathbf{y}}^{(n_1)})^\top\end{bmatrix}$;

Compute $\mathbf{f} = U_{n_3}^*\widetilde{P}_{e\mathbf{z}}/\sqrt{n_1 n_2 n_3}$ by forward FFT; $\quad \mathbf{f} = \mathbf{f}(:)$.